\documentclass[reqno]{amsart}

\usepackage{amsthm}
\usepackage{amssymb}
\usepackage{amsmath}
\usepackage{amsfonts}

\usepackage{graphicx,psfrag}

\usepackage{subcaption}

\usepackage[utf8]{inputenc}
\usepackage[english]{babel}

\usepackage[T1]{fontenc}
\usepackage{mathpazo}

\usepackage{bookmark}

\RequirePackage{amsmath}
\RequirePackage{amsfonts}

\let\:\colon

\newcommand{\cH}{\mathcal{H}}

\newcommand{\fC}{\mathfrak{C}}

\newcommand{\fM}{\mathfrak{M}}

\newcommand{\shrink}{/}

\newcommand{\deleset}{\setminus}
\newcommand{\dele}{\setminus}
\newcommand{\addeset}{+}
\newcommand{\adde}{+}
\newcommand{\delvset}{-}
\newcommand{\delv}{-}

\usepackage{xcolor}
\usepackage{hyperref}

\hypersetup{%
   pdfauthor = {F. Botler, A. Jim\'{e}nez, M. Sambinelli, Y. Wakabayashi },%
    pdftitle = {On the structure of a smallest counterexample and a new class verifying  the 2-Decomposition Conjecture},%
    pdfkeywords = {2-Decomposition Conjecture},%
    pdfproducer = {Latex with hyperref},%
    backref = page,%
    pdfcreator = {xelatex},%
    colorlinks,%
    linkcolor = {red!60!black},%
    citecolor = {green!60!black},%
    urlcolor = {blue!60!black}%
}

\usepackage{tikz,nicefrac,pifont,fp}
\usetikzlibrary{arrows,snakes,backgrounds,patterns,matrix,shapes,fit,calc,shadows,
plotmarks,arrows.meta,positioning,decorations.pathmorphing,decorations.pathreplacing,
decorations.markings,fixedpointarithmetic,hobby}

\usepackage{xifthen}
\usepackage{standalone}

\tikzset{edge/.style={line width=1.3pt, -}}
\tikzset{M edge/.style={line width=1.3pt,double distance=1.1pt,->}}
\tikzset{F1 edge/.style={line width=1.3,color=red,->}}
\tikzset{F2 edge/.style={line width=1.3,color=blue,->}}
\tikzset{E edge/.style={line width=1.3,color=black,-}}

\tikzset{red vertex/.style={circle,draw,minimum size=1mm,inner sep=0pt,outer sep=2pt,fill=red, color=red}}
\tikzset{blue vertex/.style={circle,draw,minimum size=2mm,inner sep=0pt,outer sep=4pt,fill=blue, color=blue}}
\tikzset{black vertex/.style={circle,draw,minimum size=1mm,inner sep=0pt,outer sep=2pt,fill=black, color=black}}
\tikzset{gray vertex/.style={circle,draw,minimum size=1mm,inner sep=0pt,outer sep=2pt,fill=gray, color=gray}}
\tikzset{white vertex/.style={circle,draw,minimum size=2mm,inner sep=0pt,outer sep=4pt, color=black}}

\makeatletter
\let\origsection=\section 
\def\section{\@ifstar{\origsection*}{\mysection}}
\def\mysection{\@startsection{section}{1}\z@{.7\linespacing\@plus\linespacing}{.5\linespacing}{\normalfont\scshape\centering\S}}
\makeatother

\usepackage[abbrev,msc-links,backrefs]{amsrefs}
\usepackage{doi}

\renewcommand{\PrintDOI}[1]{\doi{#1}}

\usepackage{datetime}
\usepackage{lineno}
\newcommand*\patchAmsMathEnvironmentForLineno[1]{%
\expandafter\let\csname old#1\expandafter\endcsname\csname #1\endcsname
\expandafter\let\csname oldend#1\expandafter\endcsname\csname end#1\endcsname
\renewenvironment{#1}%
{\linenomath\csname old#1\endcsname}%
{\csname oldend#1\endcsname\endlinenomath}}%
\newcommand*\patchBothAmsMathEnvironmentsForLineno[1]{%
\patchAmsMathEnvironmentForLineno{#1}%
\patchAmsMathEnvironmentForLineno{#1*}}%
\AtBeginDocument{%
\patchBothAmsMathEnvironmentsForLineno{equation}%
\patchBothAmsMathEnvironmentsForLineno{align}%
\patchBothAmsMathEnvironmentsForLineno{flalign}%
\patchBothAmsMathEnvironmentsForLineno{alignat}%
\patchBothAmsMathEnvironmentsForLineno{gather}%
\patchBothAmsMathEnvironmentsForLineno{multline}%
}

\usepackage{geometry}
\newtheorem{theorem}[equation]{Theorem}

\newtheorem{proposition}[equation]{Proposition}
\newtheorem{conjecture}[equation]{Conjecture}

\newtheorem{claim}{Claim}

\newenvironment{claimproof}[1][\proofname]
  {%
    \proof[#1]%
  }
  {%
    \endproof%
  }

\theoremstyle{definition}

\newtheorem{remark}[equation]{Remark}

\theoremstyle{case}

\numberwithin{equation}{section}

\title{On the structure of a smallest counterexample and a new class verifying  the 2-Decomposition Conjecture} 

\author{F. Botler }
\address{Programa de Engenharia de Sistemas e Computa\c{c}\~{a}o\\
Instituto Alberto Luiz Coimbra de P{\'o}s-Gradua\c{c}\~{a}o e Pesquisa em Engenharia\\
Universidade Federal do Rio de Janeiro\\
Rio de Janeiro, Brazil}  
  
\author{A. Jim\'{e}nez}
\address{Instituto de Ingenier\'ia Matem\'atica - CIMFAV\\
Facultad de Ingenier\'ia\\
Universidad de Valpara\'iso\\ Valpara\'iso, Chile} 
  
\author{M. Sambinelli }
\address{Centro de Matem{\'a}tica, Computa\c{c}\~{a}o e Cogni\c{c}\~{a}o\\
Universidade Federal do ABC\\
S\~ao Paulo, Brazil} 
  
\author{Y. Wakabayashi }
\address{Instituto de Matem\'atica e Estat\'{\i}stica\\
Universidade de S\~ao Paulo\\
S\~ao Paulo, Brazil} 
  
\thanks{%
  F. Botler is supported by CNPq (423395/2018-1) and by FAPERJ (211.305/2019);
  A. Jim\'enez is  supported by  ANID/Fondecyt Regular 1220071, MATHAMSUD 22-MATH-06, and ANID/PCI/REDES 190071;
  M. Sambinelli is supported by FAPESP (2017/23623-4) and CNPq (423833/2018-9);
  Y. Wakabayashi is supported by CNPq (306464/2016-0, 423833/2018-9) and FAPESP (2015/11937-9).
  FAPERJ and FAPESP are the Research Foundations of the states Rio de Janeiro and S\~ao Paulo, respectively.
  This study is financed in part by CAPES Finance
  Code~001 and  ANID/PCI-FAPESP 2019/13364-7.
\\
  E-mails: 
  fbotler@cos.ufrj.br (F. Botler), 
  andrea.jimenez@uv.cl (A. Jim\'enez),
  m.sambinelli@ufabc.edu.br (M. Sambinelli),
  yw@ime.usp.br (Y. Wakabayashi)
}

\newcommand{\sepcycle}[2]{\mathcal{S}_{#1, #2}}
\newcommand{\iiDCfamily}{\mathcal{S}_{2, 3}}
\newcommand{\subcubicIIDCfamily}{\mathcal{S}_{1, 3}}

\begin{document}

\maketitle

\begin{abstract} 
  The 2-Decomposition Conjecture, equivalent to the 3-Decomposition Conjecture
  stated in 2011 by Hoffmann-Ostenhof, claims that every connected graph $G$
  with vertices of degree 2 and 3, for which $G \deleset E(C)$ is disconnected
  for every cycle~$C$, admits a decomposition into a spanning tree and a
  matching. In this work we present two main results focused on developing a
  strategy to prove the 2-Decomposition Conjecture.
  One of them is a list of structural properties of a minimum counterexample
  for this conjecture. Among those properties, we prove that a minimum
  counterexample has girth at least 5 and its vertices of degree 2 are at
  distance at least 3.
  Motivated by the class of smallest counterexamples, we show that the
  2-Decomposition Conjecture holds for graphs whose vertices of degree 3 induce
  a collection of cacti in which each vertex belongs to a cycle. 
  The core of the proof of this result may possibly be used in an inductive
  proof of the 2-Decomposition Conjecture based on a parameter that relates the
  number of vertices of degree 2 and 3 in a minimum counterexample. 
\end{abstract}

\section{Introduction}

A \emph{Homeomorphically Irreducible Spanning Tree} (HIST) in a graph \(G\) is
a spanning tree of \(G\) without vertices of degree~$2$.
The problem of deciding whether a graph contains a HIST is
NP-complete~\cite{AlBeHuTh90}, even for subcubic graphs~\cite{Do92}.
This topic has been studied by many
researchers~\cite{AlBeHuTh90,ChReSh12,ChSh13,MR0351871,MR3990678} and it is
related to the topic addressed in this paper, as we shall explain.

Let $G$ be a connected cubic graph, $T$ be a spanning tree of $G$, and let $G'$
be the graph obtained from $G  \deleset  E(T)$ by removing the isolated
vertices.
Each component of $G'$ is either a path or a cycle, so every connected cubic
graph can be decomposed into a spanning tree, a collection of cycles, and a
collection of paths.
Moreover, $T$ is a HIST if and only if $G'$ is a collection of cycles.
Thus deciding whether a cubic graph $G$ contains a HIST is equivalent to
deciding whether~$G$ admits a decomposition into a spanning tree and a
collection of cycles.
Not all cubic graphs admit such a decomposition;
necessary conditions for its existence have been
shown by Hoffmann-Ostenhof, Noguchi, and Ozeki~\cite{HoNoOz18};
and the following more relaxed decomposition has been
conjectured by Hoffmann-Ostenhof~\cite{Ho11}.

\begin{conjecture}[Hoffmann-Ostenhof, 2011~\cite{Ho11}]\label{con:3dc}
  Every connected cubic graph can be decomposed into a spanning tree, a
  collection of cycles, and a (possibly empty) matching.
\end{conjecture}

Conjecture~\ref{con:3dc} is known as the \emph{$3$-Decomposition Conjecture}
(3DC, for short) and has attracted the attention of many researchers.
Although the general problem remains open, it has been verified for some
classes of cubic graphs.
Clearly, Conjecture~\ref{con:3dc} holds for every graph that contains a HIST.
Liu and Li~\cite{LiLi20} verified it for cubic traceable graphs;
and Ozeki and Ye~\cite{OzYe16} verified it for $3$-connected planar cubic
graphs and $3$-connected cubic graphs on the projective plane.
Later, Hoffmann-Ostenhof, Kaiser and Ozeki~\cite{HoKaOz18} extended the result
of Ozeki and Ye by verifying Conjecture~\ref{con:3dc} for all planar cubic
graphs.
Bachtler and Krumke~\cite{BaKr20+} verified the 3DC for a superclass of
Hamiltonian cubic graphs.
Recently Xie, Zhou and Zhou~\cite{XiZhZh20} verified Conjecture~\ref{con:3dc}
for cubic graphs containing a $2$-factor consisting of three cycles and,
independently, Hong, Liu and Yu~\cite{HoLiYu20} and Aboomahigir, Ahanjideh and
Akbari~\cite{AbAhAk20} verified the conjecture for claw-free cubic graphs.

In addition, some weaker forms of Conjecture~\ref{con:3dc} have been verified.
Akbari, Jensen, and Siggers~\cite{AkJeSi15} proved that every cubic graph can be
decomposed into a spanning forest, a collection of cycles, and a matching.
Li and Cui~\cite{LiCu14} showed that every cubic graph can be decomposed into a
spanning tree, one cycle, and a collection of paths with length at most~$2$.
Lyngsie and Merker~\cite{LyMe19} proved that every connected (not necessarily
cubic) graph can be decomposed into a spanning tree, an even graph, and a star
forest, which implies the result of Li and Cui~\cite{LiCu14} when applied to
cubic graphs.

A cycle \(C\) in a connected graph \(G\) is \emph{separating} if \(G  \deleset
E(C)\) is disconnected.
Let \(\mathcal{S}\) be the class of connected graphs in which every cycle is
separating, and let  \(\sepcycle{p}{q}\subseteq\mathcal{S}\) be the class of
graphs in \(\mathcal{S}\) with minimum degree at least $p$, and the maximum
degree at most $q$.
It is known that Conjecture~\ref{con:3dc} is equivalent to the following
conjecture, known as the \emph{$2$-Decomposition Conjecture} (2DC, for short)
-- see~\cite[Proposition~14]{HoKaOz18}.

\begin{conjecture}[Hoffmann-Ostenhof, 2016~\cite{Ho16}]\label{con:2dc}
  Every graph in \(\iiDCfamily\) can be decomposed into a spanning tree and a
  matching.
\end{conjecture}

Conjecture~\ref{con:2dc} is fairly new.
At the best of our knowledge, the only work addressing directly
Conjecture~\ref{con:2dc} is the one conducted by Hoffmann-Ostenhof, Kaiser and
Ozeki~\cite{HoKaOz18}, where the authors verify Conjecture~\ref{con:2dc} for the
planar case and show that this result implies that Conjecture~\ref{con:3dc} also
holds for planar graphs.

Throughout this paper, given a graph $G$, we denote by $V_k(G)$ the set of
vertices of~$G$ with degree~$k$.
We recall that a graph is a \emph{cactus} if it is connected and every edge is
contained in at most one cycle.
We say that a cactus~$G$ is \emph{thick} if every vertex in $G$ belongs to a
cycle.
The first contribution of this paper is the following theorem, proved in
Section~\ref{sec:maintheo}, that verifies 2DC on  graphs \(G \in \iiDCfamily\)
whose subgraph induced by \(V_3(G)\) is a collection of thick cacti.
This result was presented in a preliminary version of this work at LAGOS
2021~\cite{BoJiSaWo22+}.

\begin{theorem}\label{theo:main}
  Every graph $G \in \iiDCfamily$ for which $G  \delvset  V_2(G)$ is a
  collection of thick cacti can be decomposed into a spanning tree and a
  matching.
\end{theorem}

Let $G$ be a cubic graph and let $M$ be a perfect matching in $G$.
Note that $G - M$ is a collection of cycles, and let $G_M$ be the graph
obtained from $G$ by shrinking each cycle of \(G - M\) into a vertex.
It seems that a natural path to tackle the 3DC, and now the 2DC, is to organize
the vertices of the studied graph in a structure of cycles.
For example, it is a trivial task to verify the 3DC for Hamiltonian graphs, for
which \(G_M\) is a single vertex.
Bachtler and Krumke~\cite{BaKr20+} verified Conjecture~\ref{con:3dc} for
$3$-connected graphs containing a matching $M$ such that $G_M$ is a star;
Xie, Zhou and Zhou~\cite{XiZhZh20} verified Conjecture~\ref{con:3dc}
for  graphs containing a $2$-factor consisting of three cycles, and
here in this paper, we verify Conjecture~\ref{con:2dc} for graphs $G$ for which
$G[V_{3}]$ is a collection of cacti in which every vertex belongs to a cycle.

Another interesting and natural approach to explore these conjectures is to
study properties and forbidden structures in a minimum counterexample.
Here, we show that if \(G\) is a minimum counterexample for the 2DC, then \(G\)
is simple, 2-edge connected, has girth at least~$5$, that the distance between
any pair of vertices in \(V_2(G)\) is at least \(3\), and that its subgraph
induced by \(V_3(G)\) is connected and contains a cycle. This is precisely the
result stated in Theorem~\ref{thm:smallcount2dc} (see below), the proof of
which is given in Section~\ref{sec:counterexample}.

For a graph \(G\) in \(\sepcycle{2}{3}\) let 
\begin{equation}\label{eq:varphi}
  \varphi(G) = |V(G)| + |E(G)|.
\end{equation}
We say that a graph \(G\in\sepcycle{2}{3}\) is a \emph{counterexample} to the
2DC if \(G\) cannot be decomposed into a spanning tree and a matching.
Moreover, we say that a counterexample \(G\) is \emph{minimum} if
\(\varphi(G)\) is minimum among all counterexamples to the 2DC.

\begin{theorem}[Structure of a Minimum Counterexample to 2DC]\label{thm:smallcount2dc}
  Every minimum counterexample~$G$ to Conjecture~\ref{con:2dc} satisfies the
  following properties:
  \begin{enumerate}
    \item $G$ is simple and 2-edge-connected;
    \item the girth of $G$ is at least 5;
    \item \label{thm:smallcount2dc3}the distance between vertices in \(V_2(G)\) is at least 3;
    \item $G  \delvset  V_2(G)$ is connected;
    \item There is a cycle $C \subseteq G$ for which \(V(C)\subseteq V_3(G)\).
  \end{enumerate}
\end{theorem}

At a first glance, Theorems~\ref{theo:main} and~\ref{thm:smallcount2dc} seem
unrelated. 
In what follows we explain a connection between them.
Consider the graph parameter \(\rho\) defined over the graphs
\(G\in\sepcycle{2}{3}\) by
$$\rho(G)=|V_3(G)| - 2 |V_2(G)|.$$
By Theorem~\ref{thm:smallcount2dc}~\eqref{thm:smallcount2dc3}, if \(G\) is a
minimum counterexample to  Conjecture~\ref{con:2dc}, then \(\rho(G) \geq 0\).
The main tool in the proof of Theorem~\ref{theo:main} is
Proposition~\ref{prop:main}, which implies that graphs \(G\) for which
\(\rho(G) = 0\) admit a decomposition into a forest and a matching,
and hence may possibly be used as a base case in a proof of
Conjecture~\ref{con:2dc} by induction on~\(\rho\).

\subsection*{Notation and terminology}

The terminology used in this is work is standard and we refer the reader
to~\cite{BondyMurty2008} for missing definitions.
All graphs considered in here are finite and have no loops (but may contain
parallel edges).
Let $G = (V, E)$ be a graph.
We define  $V(G) = V$ and $E(G) = E$.
Given a vertex $u \in V(G)$, we denote the degree of $u$  by $d_G(u)$ and the neighborhood of $u$ in $G$ by $N_G(u)$ (when $G$ is clear from the context, we may simply write $d(u)$ and $N(u)$).
We denote the minimum degree of $G$ by $\delta(G)$.
We write $H \subseteq G$ to denote that $H$ is a subgraph of $G$.
Given two graphs $G$ and $H$, we write $G \cup H$ to denote the graph $(V(G) \cup V(H), E(G) \cup E(H))$.

Given a set $S \subseteq V(G)$, we denote by $G[S]$ the subgraph of $G$ induced by the vertices in $S$ and 
denote by $G \delvset S$ the graph $(V(G) \setminus S, F)$, where $F = E(G) \setminus \{uv \in E(G) \colon u \in S\}$.
When $S = \{u\}$, we may simple write $G \delvset u$ instead of $G \delvset \{u\}$.
Given a set $F \subseteq E(G)$, we denote by $G \deleset F$ the graph $(V(G), E(G) \setminus F)$, and by $G \addeset F$ the graph $(V(G), E(G) \cup F)$.
When $F = \{e\}$, we may simply write $G \deleset e$ and $G \addeset e$ instead of $G \deleset \{e\}$ and $G \addeset \{e\}$, respectively.

A \emph{path} in $G$ is a sequence $u_1u_2\cdots u_\ell$ in which $u_i \neq u_j$ for all $i \neq j$ and $u_iu_{i + 1} \in E(G)$ for $i = 1, 2, \ldots, \ell - 1$.
A \emph{cycle} in $G$ is a sequence $u_1u_2\cdots u_\ell u_1$ with $\ell \geq 2$  in which $u_i \neq u_j$ for all $i \neq j$, $u_iu_{i+1} \in E(G)$ for $i = 1, 2, \ldots \ell$, where $u_{\ell + 1} = u_1$, and $u_iu_{i+1} \neq u_ju_{j+1}$ for all $i \neq j$.
Note that this includes cycles with two edges.
Let $W = \{u_1, u_2, \ldots, u_\ell\}$ and $F = \{u_iu_{i + 1} \colon i = 1, 2, \ldots, \ell - 1\}$.
When convenient, we treat a path $P = u_1u_2 \cdots u_\ell$ and a cycle $C = u_1u_2 \cdots u_\ell u_1$ as being the subgraphs $P = (W, F)$ and $C = (W, F \cup \{ u_\ell u_1\})$, respectively.
Equivalently, a path is an acyclic connected graph with maximum degree at most~$2$ and a cycle is a connected 2-regular graph.

As usual, we say that a graph is \emph{cubic} (resp.
\emph{subcubic}) if all its vertices have degree~$3$ (resp.
at most~$3$).

\section{Proof of Theorem~\ref{theo:main}}\label{sec:maintheo}

Let $\mathcal{H} \subset \iiDCfamily$ be the set of all simple graphs \(H\) in
which \(V_2(H)\) is a stable set and every vertex in \(V_3(H)\) has precisely
one neighbor in \(V_2(H)\).
One may regard a graph \(H\) in \(\mathcal{H}\) as a graph obtained from a
cubic graph \(H'\) containing a perfect matching \(M'\) by subdividing each
edge of \(M'\) precisely once.
In particular, \(H  \delvset  V_2(H)\) is $2$-regular, as we state in the next
remark.

\begin{remark}\label{claim:comp-cycles}
  If \(H\) is a graph in \(\mathcal{H}\), then each component of \(H \delvset
  V_2(H)\) is a cycle.
\end{remark}

Let $H$ be a graph in $\mathcal{H}$.
We refer to the cycles in \(H  \delvset  V_2(H)\) as the \emph{basic cycles} of $H$.
Note that the vertices of basic cycles of $H$ define a partition of $V_3(H)$.
Let \(u \in V_2(H)\), and note that the neighbors $x$ and $y$ of $u$ belong to
basic cycles, say $C$ and $C'$, of $H$.
If \(C=C'\), then we say that the path \(P = xuy\) is a \emph{2-chord} of~\(C\).
In this case, $x$ and $y$ are called the \emph{ends} of $P$ and $u$ the
\emph{inner} vertex.
If \(C \neq C'\), then we say that \(u\) is a \emph{connector}.
In this case, we say that \(u\) \emph{joins} \(C\) and \(C'\).
Moreover, we say that two connectors are \emph{parallel} if they join the same
pair of basic cycles; and a collection of connectors~$\mathcal{C}$ of $H$ is
called \emph{simple} if it contains no pair of parallel connectors.

We define the \emph{basic cycles graph} (\emph{\emph{BC}-graph}, for short) of
$H \in \mathcal{H}$, which we denote by $\tilde{H}$, as the graph whose vertices
are the basic cycles of $H$ and in which two vertices $C, C'$ are adjacent
whenever the graph $H$ has a connector joining $C$ and $C'$.
Note that $\tilde{H}$ is connected because $H$ is connected.
Also, note that this definition ignores parallel connectors in the sense that a
pair of parallel connectors yields only one edge in~\(\tilde{H}\).
Given a collection~$\mathcal{C}$ of connectors of $H$, we define the
\emph{underlying \emph{BC}-graph} of~$\mathcal{C}$, denoted by
$\tilde{H}_{\mathcal{C}}$, as the spanning subgraph of \(\tilde{H}\) in which
two vertices \(C\) and \(C'\) are adjacent whenever there is a connector in
$\mathcal{C}$ joining $C$ and~$C'$.

We refer to a decomposition of a graph $G$ into a spanning forest $F$ and a
matching $M$ as a \emph{$2$-decomposition} of~$G$ and we denote it by the
ordered pair $(F,M)$.
Note that if a graph \(G\in\iiDCfamily\) admits a \(2\)-decomposition \((F,M)\),
then Conjecture~\ref{con:2dc} holds for \(G\) since we can complete $F$ to a
tree using edges of $M$.
Given a $2$-decomposition $\mathcal{D}= (F, M)$ of a graph $G\in \iiDCfamily$,
we say that a vertex $u \in V(G)$ is a \emph{full vertex} in $\mathcal{D}$ if
every edge of $G$ incident to  $u$ belongs to $F$.

The main result of this section (Theorem~\ref{theo:main}) is a consequence of the
following result.

\begin{proposition}\label{prop:main}
  Let $\mathcal{C}$ be a simple collection of connectors of a graph
  $H \in \mathcal{H}$.
  If $\tilde{H}_{\mathcal{C}}$ is a forest, then \(H\) admits a
  $2$-decomposition $\mathcal{D}= (F,M)$ such that each \(u\in\mathcal{C}\)
   is either a full vertex in $\mathcal{D}$ or is adjacent to a full vertex in
  \(\mathcal{D}\).
\end{proposition}

Before we present the proof of Proposition~\ref{prop:main}, we show how it
implies Theorem~\ref{theo:main}.

\begin{proof}[Proof of Theorem~\ref{theo:main}]
  The proof follows by induction on \(|E(G)|\).
  Let $G \in \iiDCfamily$ so that $G  \delvset  V_2(G)$ is a collection of thick cacti.
  The statement clearly holds for \(|E(G)|\leq 3\), so we may assume that
  \(|E(G)| \geq 4\).

  First, suppose that \(G\) is not a simple graph.
  If there are three copies of an edge, then the cycle containing any two of
  these copies is not a separating cycle, a contradiction.
  Thus, we may assume that $G$ has precisely two copies, say \(e\) and \(e'\),
  of an edge $xy$.
  If \(d(x)=2\) and \(d(y)=3\), then $y$ has degree~$1$ in $G  \delvset  V_2(G)$, and
  hence $G  \delvset  V_2(G)$ is not a collection of thick cacti, a contradiction.
  Thus, by symmetry, we may assume that \(d(x)=d(y)=3\).
  Let $u$ (resp.\ $v$) be the neighbor of $x$ (resp.\ $y$) distinct from $y$
  (resp.\ $x$).
  Since $e$ and \(e'\) form a cycle, say $C$, the graph $G \deleset E(C)$ is disconnected.
  Note that $G' = G \delvset \{x, y\} \adde \{uv\}$ is a graph in $\iiDCfamily$ and $G'  \delvset  V_2(G)$
  is a collection of thick cacti.
  Since \(|E(G')| < |E(G)|\), by the induction hypothesis \(G'\) admits a
  decomposition into a spanning tree $T'$ and a matching $M'$.
  We may assume that \(uv\in E(T')\) because \(uv\) is a cut edge of~\(G'\).
  Therefore \( \left(T' \deleset \{uv\} \addeset \{ux, e, yv\}, M' \cup \{e'\}\right)\) is a
  \(2\)-decomposition of $G$, as desired.
  Therefore, from now on, we assume that \(G\) is a simple graph.

  In what follows, we obtain a graph in $\mathcal{H}$ and a simple collection of
  connectors satisfying the hypothesis of Proposition~\ref{prop:main}.
  Let $H$ be the graph obtained from~$G$ by the following operations.
  \begin{quote}
    \begin{itemize}
      \item[(1)] replacing every path \(P_{j} = ux_1x_2\cdots x_kv\) with
            \(k\geq 1\), \(u,v\in V_3(G)\), and \mbox{\(x_i\in V_2(G)\)} for
            \(1\leq i\leq k\), by a path \(uw_{j}v\), where \(w_{j}\) is a new
            vertex (if $k=1$ we just rename $x_1$); and

      \item[(2)]subdividing once every edge of $G \delvset V_2(G)$ that does not belong to
            a cycle of \mbox{$G \delvset V_2(G)$}.
    \end{itemize}
  \end{quote}

  It is straightforward that $H \in \mathcal{H}$.
  Now, if \(\mathcal{C}\) is the set of vertices added in step (2),
  then, since $G \delvset V_2(G)$ is a collection of thick cacti, $\mathcal{C}$ is a
  simple collection of connectors of $H$ whose underlying \textrm{BC}-graph
  $\tilde{H}_{\mathcal{C}}$ is a forest.
  By Proposition~\ref{prop:main}, the graph $H$ admits a $2$-decomposition
  $\mathcal{D}=(F,M)$ in which every vertex \(u\in\mathcal{C}\) is
  a full vertex in \(\mathcal{D}\) or is adjacent to a full vertex in
  \(\mathcal{D}\).

  Now, from $(F,M)$, we obtain a $2$-decomposition $(F^{*},M^{*})$ of $G$.
  Let $xy$ be an edge in $M$.
  If $d_{H}(x) = d_{H}(y) = 3$, then we put $xy$ in $M^{*}$.
  Thus, we may assume that \(d_H(x)=2\).
  If \(x\) was added in (1), then, without loss of generality, there is a path
  \(P_{j} = ux_1x_2\cdots x_ky\) (with \(k\geq 1\), \(u,y\in V_3(G)\),
  \(x_i\in V_2(G)\) for \(1\leq i\leq k\)) that was replaced by the path
  \(uw_{j}y\), where $w_{j} = x$.
  In this case, we put \(x_ky\) in \(M^*\).
  If \(x\) was added in (2), then \(x\in\mathcal{C}\).
  Since \(x\) is not a full vertex in \(\mathcal{D}\), the vertex \(x\) must be
  adjacent to a full vertex in \(\mathcal{D}\), say \(z\), and so we put \(yz\) in
  \(M^*\).
  Since \(M\) is a transversal of the cycles of \(H\), i.e., $M$ contains an edge
  in each cycle of $H$, by construction $M^{*}$ is also a transversal of the
  cycles of \(G\), which implies that \(F^* = G  \deleset  E(M^*)\) is a forest.

  Finally, to obtain a decomposition of $G$ into a spanning tree and a matching
  from $(F^{*},M^*)$, one may find a minimal subset $S \subset M^{*}$ such that
  $F^{*} \addeset S$ is connected.
\end{proof}

For the next result, we use the following notation to refer to the
successor of a vertex in a cycle.
Given a cycle $C = w_{1}w_{2}\cdots w_{p}w_{1}$, for each $i \in \{1, \ldots,
p\}$, we denote the successor $w_{i + 1}$ (where $w_{p + 1} = w_{1}$) by
$w^{+}_{i}$.
Now we prove Proposition~\ref{prop:main} by proving the following stronger
statement.

\begin{proposition}\label{prop:main-stronger}
  Let $\mathcal{C}$ be a simple collection of connectors of a graph
  $H \in \mathcal{H}$.
  If $\tilde{H}_{\mathcal{C}}$ is a forest, then \(H\) admits a
  $2$-decomposition $\mathcal{D}= (F,M)$ such that the following holds.
  \begin{enumerate}
    \item \label{item1} $|M \cap E(C)| = 1$ for every basic cycle $C$ of $H$;
          and
    \item \label{item2} each $u \in \mathcal{C}$ is either a full vertex in
          $\mathcal{D}$ or is adjacent to a full vertex in \(\mathcal{D}\).
  \end{enumerate}
\end{proposition}
\begin{proof}
  Let \(\mathcal{C}\) and \(H \in \mathcal{H}\) be as in the statement.
  The proof follows by induction on $n = |V(H)|$.
  Since $H \in \mathcal{H}$, it follows that $V_2(H)$ is a stable set and every
  vertex in $V_3(H)$ has exactly one neighbor in \(V_2(H)\).
  Thus, we have \(|V_3(H)| = 2 |V_2(H)|\).
  If $|V_2(H)| = 1$, then $|V_3(H)| = 2$, and hence $H$ has parallel edges, a
  contradiction.
  Therefore, we may assume \(|V_2(H)|\geq 2\), which implies $|V_3(H)| \geq 4$.
  By Remark~\ref{claim:comp-cycles}, \(H  \delvset  V_2(H)\) is a collection of (basic)
  cycles.
  First, suppose that \(H\) has exactly one basic cycle, say \(C\).
  In this case, $H$ has no connectors which implies that
  \(\mathcal{C} = \emptyset\) and that every vertex in $C$ is the end of a
  2-chord.
  Let $x_1$ and $x_2$ be two adjacent vertices in $C$ which are the ends of two
  distinct 2-chords in $H$.
  Let \(x_1y_1z_1\) and \(x_2y_2z_2\) be the \(2\)-chords containing \(x_1\) and
  \(x_2\), respectively.
  For each \(y \in V_2(H)  \setminus  \{y_1,y_2\}\), let \(e_y\) be an arbitrary edge
  incident to \(y\), and let
  \begin{equation*}
    M = \{x_1x_2, y_1z_1, y_2z_2\} \cup
    \big\{e_y \colon y \in V_2(H)  \setminus  \{y_1,y_2\}\big\} \hfill
    \text{ and } \hfill F = G  \deleset  M.
  \end{equation*}
  Note that \((F,M)\) is a \(2\)-decomposition of $H$ as desired.
  Therefore, we may assume that \(H\) has at least two basic cycles.

  In what follows, we say that a basic cycle $C$ is of \emph{type~1} if no
  2-chord has both ends in~$C$; otherwise, we say $C$ is of \emph{type~2}.
  The following claim on basic cycles of type~1 arises naturally.

  \begin{claim}\label{claim:type1cutvertex}
    If $C$ is a basic cycle of type~1 in $H$, then $C$ is a cut vertex of
    $\tilde{H}$.
  \end{claim}
  \begin{claimproof}
    Since $C$ is of type~1, there are no 2-chords with ends in $C$, and since
    $H \in \cH$, the graph $H \deleset E(C)$ is disconnected.
    Let $H_{1}$ and $H_{2}$ be two distinct components of $H \deleset E(C)$.
    Note that there is no connector joining a vertex in $H_{1}$ to a vertex in
    $H_{2}$, and hence if $C_{1}$ is a basic cycle in $H_{1}$ and
    $C_{2}$ in $H_{2}$, the edge $C_{1}C_{2} \notin \tilde{H}$.
    Therefore, the vertex $C$ is a cut vertex of $\tilde{H}$.
  \end{claimproof}
  
  Since \(\tilde{H}\) is connected and
  \(\tilde{H}_{\mathcal{C}} \subseteq \tilde{H}\) is a forest, there is a
  spanning tree \(T\) of \(\tilde{H}\) such that
  \(\tilde{H}_{\mathcal{C}} \subseteq T\).
  By Claim~\ref{claim:type1cutvertex}, the leaves of $T$ are basic cycles of
  type~$2$.
  Now, let $\mathcal{C}^*$ be the collection of connectors so that
  $\tilde{H}_{\mathcal{C}^*} = T$.
  Since \(\tilde{H}_{\mathcal{C}}\subseteq T = \tilde{H}_{\mathcal{C}^*}\), we
  may assume \(\mathcal{C} \subseteq \mathcal{C}^*\).
  In what follows we prove that \(H\) admits a \(2\)-decomposition
  \(\mathcal{D} = (F,M)\) such that (a) \(|M \cap E(C)| = 1\) for every basic
  cycle \(C\) of \(H\); and (b) each \(u \in \mathcal{C}^*\) is either a full
  vertex in \(\mathcal{D}\) or is adjacent to a full vertex in \(\mathcal{D}\).
  Note that, since \(\mathcal{C} \subseteq \mathcal{C}^*\), (b) implies~(\ref{item2}),
  and hence the result follows.

  Let \(V_2(H) = \{y_1, \ldots, y_\ell\}\) and, for each \(y_i \in V_2(H)\), let
  \(x_i\) and \(z_i\) be the neighbors of \(y_i\).
  Note that~$V(H)$ is the disjoint union of the sets $\{x_i,y_i,z_i\}$ for
  $i \in \{1, \ldots, \ell\}$.
  Let \(C\) be a leaf of \(T=\tilde{H}_{\mathcal{C}^*}\), and put
  $$I = \{ y_i \in V_2(H) \colon |\{x_i, z_i \} \cap V(C)| = 1 \} \quad
  \text{and} \quad J = \{y_i \in V_2(H) \colon |\{x_i, z_i\} \cap V(C)| = 2\}.$$
  We may assume, without loss of generality, that \(x_i \in V(C)\) and \(z_i\notin V(C)\) for every \(y_i \in I\).
  Note that \(I \neq \emptyset\), otherwise either \(T\) is disconnected or \(T\) has only one vertex, namely \(C\),
  which implies that \(H\) has only one basic cycle, a contradiction.
  Thus, we may assume, without loss of generality, that \(y_1 \in \mathcal{C}^*\) and $x_1 \in V(C)$.
  Let \(C = u_1u_2 \ldots u_ku_1\), where \(u_1 = x_1\).
  Now, we split the proof into two cases depending on whether the vertex $x_1$ is adjacent to a vertex with a neighbor that belongs to $J$.

  \smallskip\noindent%
  \textbf{Case~1.}
  \emph{\(x_1\) is adjacent to a vertex with a neighbor that belongs to \(J\)
    (see Fig.~\ref{fig:examplecase1a})}.

  Suppose, without loss of generality, that \(u_2 = x_2\) and \(y_2\in J\).
  In what follows, we obtain a graph \(H' \in \mathcal{H}\) such that
  \(|V(H')| < |V(H)|\).
  Let \(H'\) be the graph obtained from \(H  \delvset  V(C)  \delvset  J\) by subdividing, for
  every \(y_i \in I\), the edge \(z_iz^+_i\), obtaining the vertex \(z'_i\), and
  adding the edge \(y_iz'_i\) (see Fig.~\ref{fig:examplecase1b}).
  Note that \(V_2(H') = V_2(H)  \setminus  J\).

  Let $B \neq C$ be a basic cycle of $H$.
  Note that if $BC \notin E(\tilde{H})$, then $B$ is a (basic) cycle in $H'$.
  On the other hand, if $BC \in E(\tilde{H})$, then $B$ is not a cycle in $H$
  but a subdivision $B'$ of $B$ is a (basic) cycle in $H'$.
  Let $\varphi(B) = B$, if $BC \notin E(\tilde{H})$, and $\varphi(B) = B'$,
  otherwise.
  Now we show that $\varphi$ is an isomorphism between $\tilde{H} - C$ and
  $\tilde{H}'$.
  It is not hard to see that, by the construction of $H'$, the function $\varphi$ is bijective.
  Now, if $XY \in E(\tilde{H} \delv C)$, then there is a connector, say~$y$, joining
  the basic cycles $X$ and $Y$ in $H$.
  The only connectors affected by the construction of $H'$ are those that
  contain an end in $C$, and since $X \neq C$ and $Y \neq C$, it follows that
  $y$ is a connector joining $\varphi(X)$ to $\varphi(Y)$ in $H'$.
  Thus $\varphi(X)\varphi(Y) \in E(\tilde{H'})$.
  Now, suppose that $\varphi(X)\varphi(Y) \in E(\tilde{H'})$.
  By the construction of $H'$, no connector is created.
  This implies that the set of connectors of $H'$ is a subset of the set of connectors of $H$.
  Thus, if $\varphi(X)\varphi(Y) \in E(\tilde{H'})$, then there is a
  connector, say $y$, joining the basic cycles $\varphi(X)$ and
  $\varphi(Y)$.
  Since $y$ is also a connector in $H$, $X \neq C$, and $Y \neq C$, it follows
  that $XY \in E(\tilde{H} \delv C)$.
  Therefore, $\tilde{H} \delv C$ and $\tilde{H'}$ are isomorphic.

  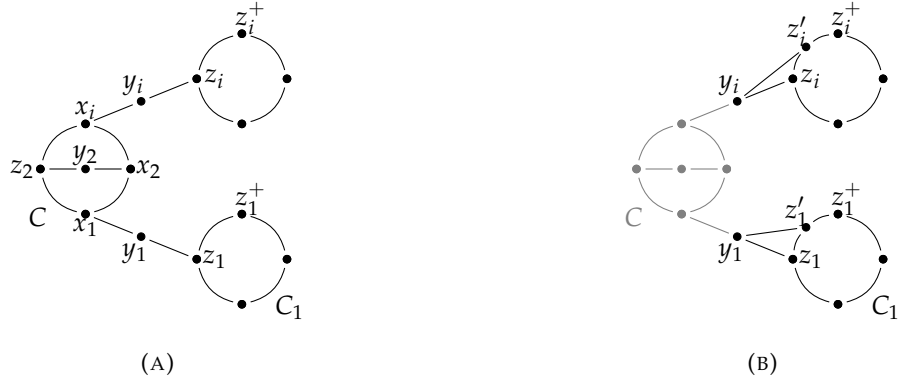
\begin{figure}
    \centering
    \begin{subfigure}{.5\textwidth}
      \centering
      \begin{tikzpicture}[scale = .6]
      
      \node (center1) []	at (0,0) {};
      \node (center2) []	at ($(center1) + (30:4)$) {};
      \node (center3) []	at ($(center1) + (-30:4)$) {};

    \foreach \j in {1,...,4}
        {
          \node(1 \j) [black vertex] at ($(center1)+(\j*90:1)$) {};
          \node(aux) [black vertex] at ($(center1)+(\j*90+90:1)$) {};
          \draw (1 \j) to [bend right=30] (aux);
        }
        
    \foreach \j in {1,...,4}
        {
          \node(2 \j) [black vertex] at ($(center2)+(\j*90:1)$) {};
          \node(aux) [black vertex] at ($(center2)+(\j*90+90:1)$) {};
          \draw (2 \j) to [bend right=30] (aux);
        }
        
    \foreach \j in {1,...,4}
        {
          \node(3 \j) [black vertex] at ($(center3)+(\j*90:1)$) {};
          \node(aux) [black vertex] at ($(center3)+(\j*90+90:1)$) {};
          \draw (3 \j) to [bend right=30] (aux);
        }


      \node(in1) [black vertex] at ($.5*(1 4) + .5*(1 2)$) {};
      \draw[] 	(1 4) -- (in1);
      \draw[]		(in1) -- (1 2);


      \node(mid12) [black vertex] at ($.5*(1 1)	+.5*(2 2)$) {};
      \draw[] 	(1 1) -- (mid12);
      \draw[] 	(mid12) -- (2 2);

      \node(mid13) [black vertex] at ($.5*(1 3)	+.5*(3 2)$) {};
      \draw[] 	(1 3) -- (mid13);
      \draw[]		(mid13) -- (3 2);
      
      \node () [] at ($(1 3) + (-90:.3)$) {$x_1$};
      \node () [] at ($(1 4) + (0:.4)$) {$x_2$};
      \node () [] at ($(1 1) + (90:.3)$) {$x_i$};
      \node () [] at ($(1 2) + (180:.4)$) {$z_2$};
      \node () [] at ($(in1) + (90:.3)$) {$y_2$};
      
      \node () [] at ($(mid12) + (115:.4)$) {$y_i$};
      \node () [] at ($(2 2) + (0:.4)$) {$z_i$};
      \node () [] at ($(2 1) + (60:.4)$) {$z_i^+$};
      
      \node () [] at ($(mid13) + (250:.4)$) {$y_1$};
      \node () [] at ($(3 2) + (0:.4)$) {$z_1$};
      \node () [] at ($(3 1) + (60:.4)$) {$z_1^+$};
      
      \node () [] at ($(center1) + (225:1.5)$) {$C$};
      \node () [] at ($(center3) + (-45:1.5)$) {$C_1$};
\end{tikzpicture}

      \caption{}
      \label{fig:examplecase1a}
    \end{subfigure}%
    \begin{subfigure}{.5\textwidth}
      \centering
      \begin{tikzpicture}[scale = .6]
	
	\node (center1) []	at (0,0) {};
	\node (center2) []	at ($(center1) + (30:4)$) {};
	\node (center3) []	at ($(center1) + (-30:4)$) {};

\foreach \j in {1,...,4}
		{
			\node(1 \j) [gray vertex] at ($(center1)+(\j*90:1)$) {};
			\node(aux) [gray vertex] at ($(center1)+(\j*90+90:1)$) {};
			\draw[color=gray] (1 \j) to [bend right=30] (aux);
		}
		
\foreach \j in {1,...,4}
		{
			\node(2 \j) [black vertex] at ($(center2)+(\j*90:1)$) {};
			\node(aux) [black vertex] at ($(center2)+(\j*90+90:1)$) {};
			\ifthenelse{\NOT \j = 1}{
				\draw (2 \j) to [bend right=30] (aux);
			}{}
		}		
		\node(2 a) [black vertex] at ($(center2)+(1.5*90:1)$) {};
		\draw (2 1) to [bend right=15] (2 a);
		\draw (2 a) to [bend right=15] (2 2);
		
\foreach \j in {1,...,4}
		{
			\node(3 \j) [black vertex] at ($(center3)+(\j*90:1)$) {};
			\node(aux) [black vertex] at ($(center3)+(\j*90+90:1)$) {};
			\ifthenelse{\NOT \j = 1}{
				\draw (3 \j) to [bend right=30] (aux);
			}{}
		}
		
		\node(3 a) [black vertex] at ($(center3)+(1.5*90:1)$) {};
		\draw (3 1) to [bend right=15] (3 a);
		\draw (3 a) to [bend right=15] (3 2);

	\node(in1) [gray vertex] at ($.5*(1 4) + .5*(1 2)$) {};
	\draw[color=gray] 	(1 4) -- (in1);
	\draw[color=gray]	(in1) -- (1 2);


	\node(mid12) [black vertex] at ($.5*(1 1)	+.5*(2 2)$) {};
	\draw[color=gray]	(1 1) -- (mid12);
	\draw[] 			(mid12) -- (2 2);
	\draw[]				(2 a) -- (mid12);

	\node(mid13) [black vertex] at ($.5*(1 3)	+.5*(3 2)$) {};
	\draw[color=gray]	(1 3) -- (mid13);
	\draw[]			 	(mid13) -- (3 2);
	\draw[]				(3 a) -- (mid13);
	
	
	\node () [] at ($(mid12) + (115:.4)$) {$y_i$};
	\node () [] at ($(2 2) + (0:.4)$) {$z_i$};
	\node () [] at ($(2 1) + (60:.4)$) {$z_i^+$};
	
	\node () [] at ($(mid13) + (250:.4)$) {$y_1$};
	\node () [] at ($(3 2) + (0:.4)$) {$z_1$};
	\node () [] at ($(3 1) + (60:.4)$) {$z_1^+$};
	
	\node () [] at ($(2 a) + (120:.4)$) {$z_i'$};
	\node () [] at ($(3 a) + (120:.4)$) {$z_1'$};
	
	\node () [] at ($(center1) + (225:1.5)$) {\color{gray}$C$};
	\node () [] at ($(center3) + (-45:1.5)$) {$C_1$};

\end{tikzpicture}
      \caption{}
      \label{fig:examplecase1b}
    \end{subfigure}%
    \caption{Reduction from a graph $H$ (\ref{fig:examplecase1a}) to the graph
      $H'$ \eqref{fig:examplecase1b} in Case 1.
      In~\eqref{fig:examplecase1b}, we use gray to indicate the elements from
      $H$ that we removed to create $H'$.}
    \label{fig:examplecase1}
  \end{figure}

  Let $W = V(H) \setminus (I \cup J \cup V(C))$, and note that
  $V(H) = W \cup I \cup J \cup V(C)$.
  Moreover, note that $V(H') = W \cup I \cup \{z'_{i} \colon y_{i} \in I\}$.
  Since \(J \neq \emptyset\), it follows that $|V(H')|< |V(H)|$.
  Now we show that \(H' \in \mathcal{H}\).
  First, since $C$ is a leaf of the spanning tree $T$ of $\tilde{H}$, the graph
  $T \delv C$ is connected, and since $\tilde{H'}$ is isomorphic to $\tilde{H} \delv C$,
  we have that $T \delv C$ is a spanning tree of $\tilde{H'}$, and hence $H'$ is
  connected.
  Also, by construction, $H'$ is simple, each vertex in \(V_3(H')\) has exactly
  one neighbor in \(V_2(H')\), and \(V_2(H')\) is a stable set.
  It remains to prove that every cycle in \(H'\) is a separating cycle.
  First, note that every basic cycle $C'$ adjacent to $C$ in $\tilde{H}$ yields
  a basic cycle of type~2 in $\tilde{H}'$.
  Moreover, note that a cycle containing a vertex with degree~$2$ or a basic
  cycle of type~2 is a separating cycle.
  Thus, we can focus on the basic cycles of type~1 in $H'$.
  Let $C'$ be such a cycle.
  By Claim~\ref{claim:type1cutvertex} the cycle $C'$ is a cut vertex of
  $\tilde{H}$.
  By the construction of $H'$, the graph $\tilde{H} \delv C$ is isomorphic to
  $\tilde{H'}$.
  Thus, $C'$ is a cut vertex of $\tilde{H'}$, which implies that $C'$ is a
  separating cycle of $H'$.
  Therefore, we conclude that $H' \in \mathcal{H}$.

  By the induction hypothesis, the graph \(H'\) has a \(2\)-decomposition
  \((F', M')\) satisfying~(\ref{item1}) and~(\ref{item2}) with respect to the
  simple collection of connectors $\mathcal{C}^* \setminus \{y_1\}$.
  We now describe in four steps how to obtain the desired \(2\)-decomposition
  \((F, M)\) of $H$ from \((F', M')\) (see Fig.~\ref{fig:decompcase1}):

  \begin{itemize}
    \item[(a)] We put $x_1x_2$ in $M$ and all the edges of $C  \deleset  \{x_1x_2\}$ in $F$.
          We put $y_2x_2$ in $F$, $y_2z_2$ in $M$ and, for each
          $y_i \in J \setminus \{y_2\}$, we put $x_iy_i$ and $y_iz_i$ in distinct elements of $(F, M)$.
    \item[(b)] We put $x_1y_1$ and $y_1z_1$ in $F$.
          In addition, for each $y_i \in I \setminus \{y_1\}$, we put $x_iy_i$
          in $M$ and $y_iz_i$ in $F$.
    \item[(c)] We put each
          edge $$e \in E(G)  \deleset  \Big( E(C) \cup \{y_ix_i, y_iz_i \colon y_i \in I \cup J\} \cup \{z_iz^+_i \colon y_i \in I\}\Big) $$
          in $F$ if $e \in F'$.
          Otherwise, we put $e$ in $M$.
    \item[(d)] Finally, for each $y_i \in I$, we put
          $z_iz_i^+ \in E(G)$ in $M$ if $z'_iz^+_i \in M'$.
          Otherwise, we put $z_iz_i^+$ in $F$.
  \end{itemize}

  We claim that in  steps (a)-(d), each edge in $E(H)$ has been put
  either in $M$ or $F$. Indeed, in steps~(a) and~(b) we cover the edges in 
  $E(C) \cup \{y_iz_i, y_ix_i \colon y_i \in I \cup J\}$, while the edges in $\{z_iz^+_i \colon y_i \in I \}$ are covered in step
  (d), and all
  the remaining edges are covered in step~(c).
  The following claim is useful.

  \begin{claim}\label{claim:notinM}
    Edge $z_iz'_i \in E(F')$ for all $y_i \in I$.
  \end{claim}
  \begin{claimproof}
    Let $B$ be the basic cycle of $H'$ that contains $z_iz'_i$.
    If $z_iz'_i$ belongs to $M'$, then due to~(\ref{item1}) all the edges in
    $E(B)  \deleset  \{z_iz'_i\}$ belong to $F'$.
    Hence, the cycle $\{z_iy_i, y_iz'_i\} \cup ( E(B)  \deleset  \{z_iz'_i\})$ is contained in 
    $F'$, a contradiction.
 \end{claimproof}

 It is straightforward from the assignments in  steps (a)-(d) that $F$ is a
 forest and $M$ is a matching (for the edges that were subdivided,
 item~(\ref{item1}) with respect to $M'$ and Claim~\ref{claim:notinM} ensure
 that $M$ is a matching and that $F$ has no cycles).
 We now check that items~(\ref{item1}) and~(\ref{item2}) hold for $(F,M)$.
 Due to step~(a), we have that $|M \cap E(C)|=1$ and due to steps~(c)-(d), we
 have that $|M \cap E(C')|=1$ for every other basic cycle $C'$ of $H$ with
 $C'\neq C$.
 Hence, (\ref{item1}) holds.
 Due to step~(b), we have that $y_1$ is a full vertex in $(F,M)$.
 Let $y \in \mathcal{C}^* \setminus \{y_1\}$.
 Since $C$ is a leaf of $T = \tilde{H}_{\mathcal{C}^*}$, it follows that
 $y \notin I$.
 Thus, if $y$ is full in $(F',M')$, then due to the step~(c), it is also full in
 $(F,M)$.
 So, suppose that $y$ is not full in $(F',M')$ and hence, it has a neighbor $x$
 which is a full vertex in $(F',M')$.
 Let $C'$ be the basic cycle of $H$ that contains $x$.
 If $x \neq z^+_i$ for each $y_i \in I \setminus \{y_1\}$, then due to step~(c), we have
 that $x$ is a full vertex in $(F,M)$ as well.
 Suppose that $x = z^+_i$ for some $i \in I \setminus \{y_1\}$.
 Since $x$ is a full vertex in $(F',M')$, we have that $xy, xz'_i$ and $xx'$
 belong to $F'$, where $x'$ is the neighbor of $x$ distinct of $z'_i$ in the
 basic cycle $C'$.
 Due to step~(c), $xy$ and $xx'$ belong to $F$, and due to step~(d) the edge
 $z_iz^+_i$ belongs to $F$ (since $z'_ix \in E(F')$).
 Therefore item~(\ref{item2}) holds.
 This finishes the proof of Case~1.

 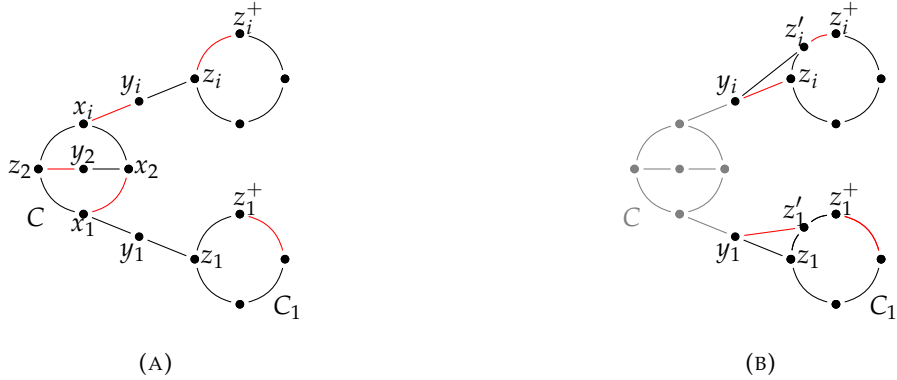
\begin{figure}
   \centering
   \begin{subfigure}{.5\textwidth}
     \centering 
     \begin{tikzpicture}[scale = .6]
	
	\node (center1) []	at (0,0) {};
	\node (center2) []	at ($(center1) + (30:4)$) {};
	\node (center3) []	at ($(center1) + (-30:4)$) {};

\foreach \j in {1,...,4}
		{
			\node(1 \j) [black vertex] at ($(center1)+(\j*90:1)$) {};
			\node(aux) [black vertex] at ($(center1)+(\j*90+90:1)$) {};
			\ifthenelse{\NOT \j = 3}{
				\draw (1 \j) to [bend right=30] (aux);
			}{
				\draw[color=red] (1 \j) to [bend right=30] (aux);
			}
		}

\foreach \j in {1,...,4}
		{
			\node(2 \j) [black vertex] at ($(center2)+(\j*90:1)$) {};
			\node(aux) [black vertex] at ($(center2)+(\j*90+90:1)$) {};
			\ifthenelse{\NOT \j = 1}{
				\draw (2 \j) to [bend right=30] (aux);
			}{
				\draw[color=red] (2 \j) to [bend right=30] (aux);
			}
		}
		
\foreach \j in {1,...,4}
		{
			\node(3 \j) [black vertex] at ($(center3)+(\j*90:1)$) {};
			\node(aux) [black vertex] at ($(center3)+(\j*90+90:1)$) {};
			\ifthenelse{\NOT \j = 4}{
				\draw (3 \j) to [bend right=30] (aux);
			}{
				\draw[color=red] (3 \j) to [bend right=30] (aux);
			}
		}


	\node(in1) [black vertex] at ($.5*(1 4) + .5*(1 2)$) {};
	\draw 	(1 4) -- (in1);
	\draw[color=red]		(in1) -- (1 2);


	\node(mid12) [black vertex] at ($.5*(1 1)	+.5*(2 2)$) {};
	\draw[color=red] 	(1 1) -- (mid12);
	\draw[] 	(mid12) -- (2 2);

	\node(mid13) [black vertex] at ($.5*(1 3)	+.5*(3 2)$) {};
	\draw[] 	(1 3) -- (mid13);
	\draw[]		(mid13) -- (3 2);
	
	\node () [] at ($(1 3) + (-90:.3)$) {$x_1$};
	\node () [] at ($(1 4) + (0:.4)$) {$x_2$};
	\node () [] at ($(1 1) + (90:.3)$) {$x_i$};
	\node () [] at ($(1 2) + (180:.4)$) {$z_2$};
	\node () [] at ($(in1) + (90:.3)$) {$y_2$};
	
	\node () [] at ($(mid12) + (115:.4)$) {$y_i$};
	\node () [] at ($(2 2) + (0:.4)$) {$z_i$};
	\node () [] at ($(2 1) + (60:.4)$) {$z_i^+$};
	
	\node () [] at ($(mid13) + (250:.4)$) {$y_1$};
	\node () [] at ($(3 2) + (0:.4)$) {$z_1$};
	\node () [] at ($(3 1) + (60:.4)$) {$z_1^+$};
	
	\node () [] at ($(center1) + (225:1.5)$) {$C$};
	\node () [] at ($(center3) + (-45:1.5)$) {$C_1$};
\end{tikzpicture}

     \caption{}
     \label{fig:examplecase1c}
   \end{subfigure}%
   \begin{subfigure}{.5\textwidth}
     \centering 
     \begin{tikzpicture}[scale = .6]

	\node (center1) []	at (0,0) {};
	\node (center2) []	at ($(center1) + (30:4)$) {};
	\node (center3) []	at ($(center1) + (-30:4)$) {};

\foreach \j in {1,...,4}
		{
			\node(1 \j) [gray vertex] at ($(center1)+(\j*90:1)$) {};
			\node(aux) [gray vertex] at ($(center1)+(\j*90+90:1)$) {};
			\draw[color=gray] (1 \j) to [bend right=30] (aux);
		}

\foreach \j in {1,...,4}
		{
			\node(2 \j) [black vertex] at ($(center2)+(\j*90:1)$) {};
		}

\foreach \j in {1,...,4}
		{
			\node(aux) [black vertex] at ($(center2)+(\j*90+90:1)$) {};
			\ifthenelse{\NOT \j = 1}{
				\draw[] (2 \j) to [bend right=30] (aux);
			}{
				\node(2 a) [black vertex] at ($(center2)+(1.5*90:1)$) {};
				\draw[color=red] (2 1) to [bend right=15] (2 a);
				\draw (2 a) to [bend right=15] (2 2);
			}
		}

\foreach \j in {1,...,4}
		{
			\node(3 \j) [black vertex] at ($(center3)+(\j*90:1)$) {};
		}

\foreach \j in {1,...,4}
		{
			\node(3 \j) [black vertex] at ($(center3)+(\j*90:1)$) {};
			\node(aux) [black vertex] at ($(center3)+(\j*90+90:1)$) {};
			\ifthenelse{\NOT \j = 1 \AND \NOT \j = 4}{
				\draw (3 \j) to [bend right=30] (aux);
			}{
				\draw[color=red] (3 4) to [bend right=30] (3 1);
				\node(3 a) [black vertex] at ($(center3)+(1.5*90:1)$) {};
				\draw (3 1) to [bend right=15] (3 a);
				\draw (3 a) to [bend right=15] (3 2);
			}
		}


	\node(in1) [gray vertex] at ($.5*(1 4) + .5*(1 2)$) {};
	\draw[color=gray] 	(1 4) -- (in1);
	\draw[color=gray]	(in1) -- (1 2);

Bridges

	\node(mid12) [black vertex] at ($.5*(1 1)	+.5*(2 2)$) {};
	\draw[color=gray]	(1 1) -- (mid12);
	\draw[color=red] 			(mid12) -- (2 2);
	\draw				(2 a) -- (mid12);

	\node(mid13) [black vertex] at ($.5*(1 3)	+.5*(3 2)$) {};
	\draw[color=gray]	(1 3) -- (mid13);
	\draw			 	(mid13) -- (3 2);
	\draw[color=red]				(3 a) -- (mid13);


	\node () [] at ($(mid12) + (115:.4)$) {$y_i$};
	\node () [] at ($(2 2) + (0:.4)$) {$z_i$};
	\node () [] at ($(2 1) + (60:.4)$) {$z_i^+$};

	\node () [] at ($(mid13) + (250:.4)$) {$y_1$};
	\node () [] at ($(3 2) + (0:.4)$) {$z_1$};
	\node () [] at ($(3 1) + (60:.4)$) {$z_1^+$};

	\node () [] at ($(2 a) + (120:.4)$) {$z_i'$};
	\node () [] at ($(3 a) + (120:.4)$) {$z_1'$};

	\node () [] at ($(center1) + (225:1.5)$) {\color{gray} $C$};
	\node () [] at ($(center3) + (-45:1.5)$) {$C_1$};
\end{tikzpicture}
     \caption{}
     \label{fig:examplecase1d}
   \end{subfigure}%

   \caption{Reduction from a graph $H$ \eqref{fig:examplecase1c} to the graph $H'$
     \eqref{fig:examplecase1d} in Case 1.
     In both figures, for a \(2\)-decomposition \((F,M)\), the edges in \(F\)
     (resp.\ \(M\)) are colored black (resp.\ red).
     In~\eqref{fig:examplecase1d}, we use gray to indicate the elements from $H$
     that we removed to create $H'$.}%
   \label{fig:decompcase1}

 \end{figure}

 \medskip
 \noindent
 \textbf{Case~2.}
 \emph{\(x_1\) is not adjacent to a vertex with a neighbor in \(J\).}
  
 In this case, both neighbors of \(x_1\) in \(C\), namely \(u_2\) and \(u_k\),
 are adjacent to a vertex in \(I\).
 Let $C_1$ be the basic cycle that contains $z_1$ and let $\ell$ be the smallest
 $i \in\{1, \ldots, k\}$ for which $u_{i+1}$ is an end of a \(2\)-chord.
 We may assume, without loss of generality, that $u_{\ell}=x_2$ and
 $u_{\ell +1}=x_3$.
 Let $C_2$ be the basic cycle that contains $z_2$ (possibly $C_{1} = C_{2}$).

 In what follows, analogously to Case 1, we obtain a graph
 \(H' \in \mathcal{H}\) such that \(|V(H')|<|V(H)|\).
 Let \(H'\) be the graph obtained from \(H  \delvset  V(C)  \delvset  J\) by (i) identifying the
 vertices \(y_1\) and \(y_2\) into a new vertex \(y\), and (ii) for every
 \(y_i \in I \setminus \{y_1,y_2\}\), subdividing the edge \(z_iz^+_i\), obtaining the
 vertex \(z'_i\), and adding the edge \(y_iz'_i\) (see
 Fig.~\ref{fig:examplecase2}).
 Now, note that
 \(V_2(H') = \{y\} \cup \big(V_2(H) \setminus \left(J \cup \{y_1, y_2\}\right)\big)\).

 Let $B \neq C$ be a basic cycle of $H$.
 Note that if $BC \notin E(\tilde{H})$, then $B$ is a (basic) cycle in $H'$.
 On the other hand, if $BC \in E(\tilde{H})$, then $B$ is not a cycle in $H'$,
 but one subdivision $B'$ of $B$ is.
 Let $\varphi(B) = B$, if $BC \notin E(\tilde{H})$, and $\varphi(B) = B'$,
 otherwise.
 Now we show that, if $C_{1} = C_{2}$, then $\varphi$ is an isomorphism between
 $\tilde{H'}$ and $\tilde{H} \delv C$; otherwise, we show that $\varphi$ is an
 isomorphism between $\tilde{H'}$ and $\tilde{H} \delv C \adde C_{1}C_{2}$ (here,
 we only add the edge $C_{1}C_{2}$ to $\tilde{H} \delv C$ if this action results in
 a simple graph).
 It is not hard to check that, by the construction of $H'$, 
 the function $\varphi$ is bijective.
 If $XY \in E(\tilde{H} \delv C)$, then there is a connector, say~$y'$, joining the
 basic cycles $X$ and $Y$ in $H$.
 The only connectors affected by the construction of $H'$ are those that contain
 an end in $C$ and, since $X \neq C$ and $Y \neq C$, it follows that $y'$ is a
 connector joining $\varphi(X)$ to $\varphi(Y)$ in $H'$.
 If $C_{1} \neq C_{2}$, then, by the construction of $H'$, the vertex $y$ is a
 connector in $H'$ joining the basic cycles $\varphi(C_{1})$ and
 $\varphi(C_{2})$ and, as result,
 $\varphi(C_{1})\varphi(C_{2}) \in E(\tilde{H'})$.
 Now, suppose that $\varphi(X)\varphi(Y) \in E(\tilde{H'})$, and hence there
 exists a connector $y'$ in $H'$ joining $\varphi(X)$ and $\varphi(Y)$.
 By the construction of $H'$, the vertex $y$ is the only connector that we can
 create, so the set of connectors of $H'$ distinct from $y$ is a subset of
 the set of connectors of $H$.
 If $|\{X, Y\} \cap \{C_{1}, C_{2}\}| < 2$, then $y' \neq y$, and hence $y'$ is
 a connector in $H$ joining $X$ and $Y$, and hence $XY \in E(\tilde{H} \delv C)$.
 Now, $|\{X, Y\} \cap \{C_{1}, C_{2}\}| = 2$, then, by construction,
 $XY \in \tilde{H} \delv C \adde C_{1}C_{2}$.
 Therefore, $\varphi$ is an isomorphism between $\tilde{H'}$ and
 $\tilde{H} \delv C$, if $C_{1} = C_{2}$, or between $\tilde{H'}$ and
 $\tilde{H} \delv C \adde C_{1}C_{2}$, otherwise.

\begin{figure}
\centering
\begin{subfigure}{.5\textwidth}
    \centering
    \begin{tikzpicture}[scale = .6]
	
	\node (center1) []	at (0,0) {};
	\node (center2) []	at ($(center1) + (36:4)$) {};
	\node (center3) []	at ($(center1) + (-36:4)$) {};
	\node (center4) []	at ($(center1) + (-72-36:4)$) {};

\foreach \j in {1,...,5}
		{
			\node(1 \j) [black vertex] at ($(center1)+(\j*72+36:1)$) {};
			\node(aux) [black vertex] at ($(center1)+(\j*72+72+36:1)$) {};
			\draw (1 \j) to [bend right=18] (aux);

		}
		
\foreach \j in {1,...,4}
		{
			\node(2 \j) [black vertex] at ($(center2)+(\j*90+36:1)$) {};
			\node(aux) [black vertex] at ($(center2)+(\j*90+90+36:1)$) {};
			\draw (2 \j) to [bend right=30] (aux);
		}
		
\foreach \j in {1,...,4}
		{
			\node(3 \j) [black vertex] at ($(center3)+(\j*90-36:1)$) {};
			\node(aux) [black vertex] at ($(center3)+(\j*90+90-36:1)$) {};
			\draw (3 \j) to [bend right=30] (aux);
		}
		
\foreach \j in {1,...,4}
		{
			\node(4 \j) [black vertex] at ($(center4)+(\j*90-18:1)$) {};
			\node(aux) [black vertex] at ($(center4)+(\j*90+90-18:1)$) {};
			\draw (4 \j) to [bend right=30] (aux);
		}


	\node(in1) [black vertex] at ($.5*(1 1) + .5*(1 2) + (-45:.4)$) {};
	\draw[] 	(1 1) -- (in1);
	\draw[]		(in1) -- (1 2);


	\node(mid12) [black vertex] at ($.5*(1 5)	+.5*(2 2)$) {};
	\draw[] 	(1 5) -- (mid12);
	\draw[] 	(mid12) -- (2 2);

	\node(mid13) [black vertex] at ($.5*(1 4)	+.5*(3 2)$) {};
	\draw[] 	(1 4) -- (mid13);
	\draw[]		(mid13) -- (3 2);
	
	\node(mid14) [black vertex] at ($.5*(1 3)	+.5*(4 1)$) {};
	\draw[] 	(1 3) -- (mid14);
	\draw[]		(mid14) -- (4 1);

	\node () [] at ($(1 3) + (60:.3)$) {$u_k$};
	\node () [] at ($(1 4) + (-90:.4)$) {$x_1$};
	\node () [] at ($(1 5) + (90:.4)$) {$x_2$};
	\node () [] at ($(1 1) + (100:.4)$) {$x_3$};
	\node () [] at ($(1 2) + (180:.4)$) {$z_3$};
	\node () [] at ($(in1) + (-30:.4)$) {$y_3$};
	
	\node () [] at ($(mid12) + (115:.4)$) {$y_2$};
	\node () [] at ($(2 2) + (45:.4)$) {$z_2$};
	
	\node () [] at ($(mid13) + (250:.4)$) {$y_1$};
	\node () [] at ($(3 2) + (-45:.4)$) {$z_1$};

	\node () [] at ($(center1) + (225:1.5)$) {$C$};
	\node () [] at ($(center3) + (-45:1.5)$) {$C_1$};
	\node () [] at ($(center2) + (-45:1.5)$) {$C_2$};
\end{tikzpicture}

    \caption{}
    \label{fig:examplecase2a}
\end{subfigure}%
\begin{subfigure}{.5\textwidth}
    \centering
\begin{tikzpicture}[scale = .6]
	
	\node (center1) []	at (0,0) {};
	\node (center2) []	at ($(center1) + (36:4)$) {};
	\node (center3) []	at ($(center1) + (-36:4)$) {};
	\node (center4) []	at ($(center1) + (-72-36:4)$) {};

\foreach \j in {1,...,5}
		{
			\node(1 \j) [gray vertex] at ($(center1)+(\j*72+36:1)$) {};
			\node(aux) [gray vertex] at ($(center1)+(\j*72+72+36:1)$) {};
			\draw[color=gray] (1 \j) to [bend right=18] (aux);
		}
		
\foreach \j in {1,...,4}
		{
			\node(2 \j) [black vertex] at ($(center2)+(\j*90+36:1)$) {};
			\node(aux) [black vertex] at ($(center2)+(\j*90+90+36:1)$) {};
			\draw (2 \j) to [bend right=30] (aux);
		}
		
\foreach \j in {1,...,4}
		{
			\node(3 \j) [black vertex] at ($(center3)+(\j*90-36:1)$) {};
			\node(aux) [black vertex] at ($(center3)+(\j*90+90-36:1)$) {};
			\draw (3 \j) to [bend right=30] (aux);
		}
		
\foreach \j in {1,...,5}
		{
			\node(4 \j) [black vertex] at ($(center4)+(\j*72:1)$) {};
			\node(aux) [black vertex] at ($(center4)+(\j*72+72:1)$) {};
			\draw (4 \j) to [bend right=18] (aux);
		}


	\node(in1) [gray vertex] at ($.5*(1 1) + .5*(1 2) + (-45:.4)$) {};
	\draw[color=gray] 	(1 1) -- (in1);
	\draw[color=gray]		(in1) -- (1 2);


	\node(mid12) [gray vertex] at ($.5*(1 5)	+.5*(2 2)$) {};
	\draw[color=gray] 	(1 5) -- (mid12);
	\draw[color=gray] 	(mid12) -- (2 2);

	\node(mid13) [gray vertex] at ($.5*(1 4)	+.5*(3 2)$) {};
	\draw[color=gray] 	(1 4) -- (mid13);
	\draw[color=gray]		(mid13) -- (3 2);
	
	\node(mid14) [black vertex] at ($.5*(1 3)	+.5*(4 1)$) {};
	\draw[color=gray] 	(1 3) -- (mid14);
	\draw[]		(mid14) -- (4 1);
	\draw[]		(mid14) -- (4 5);
	
	\node(mid23) [black vertex] at ($.5*(mid12)	+.5*(mid13)$) {};
	\draw[] 	(2 2) -- (mid23);
	\draw[]		(mid23) -- (3 2);
	
	
	\node () [] at ($(2 2) + (45:.4)$) {$z_2$};
	
	\node () [] at ($(3 2) + (-45:.4)$) {$z_1$};

	\node () [] at ($(mid23) + (0:.4)$) {$y$};
	
	\node () [] at ($(center1) + (225:1.5)$) {\color{gray} $C$};
	\node () [] at ($(center3) + (-45:1.5)$) {$C_1$};
	\node () [] at ($(center2) + (-45:1.5)$) {$C_2$};
\end{tikzpicture}
    \caption{}
    \label{fig:examplecase2b}
\end{subfigure}%

\caption{Reduction from a graph $H$ \eqref{fig:examplecase2a} to the graph
    $H'$ \eqref{fig:examplecase2b} in Case 2.
    In~\eqref{fig:examplecase2b}, we use gray to indicate the elements from
    $H$ that we removed to create $H'$.}%
\label{fig:examplecase2}
\end{figure}
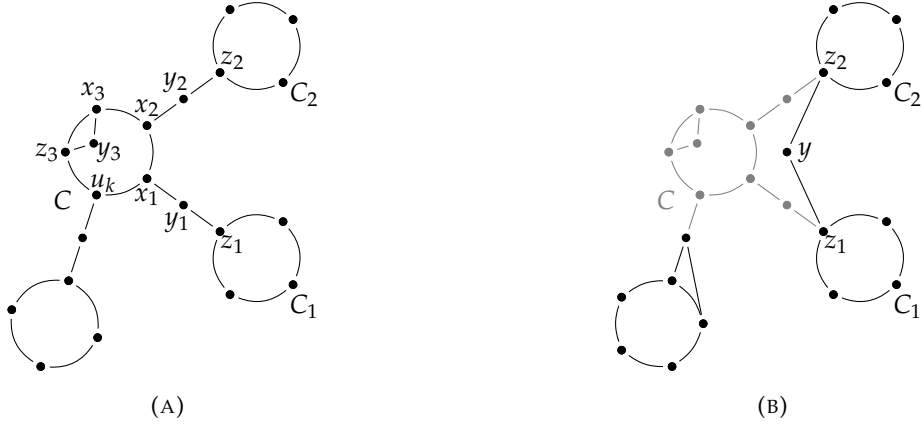

Let $W = V(H) \setminus (I \cup J \cup V(C))$, and hence
$V(H) = W \cup I \cup J \cup V(C)$.
Moreover, note that
$V(H') = W \cup (I \setminus \{y_{1}, y_{2}\}) \cup \{z'_{i} \colon y_{i} \in I \setminus \{y_{1}, y_{2}\}\} \cup \{y\}$.
It follows that $|V(H')|< |V(H)|$.
Now we claim that \(H' \in \mathcal{H}\).
First, since $C$ is a leaf of a spanning tree $T$ of $\tilde{H}$, and either
$\tilde{H} \delv C$ or $\tilde{H} \delv C \adde C_{1}C_{2}$ is isomorphic to
$\tilde{H'}$, it follows that $T \delv C$ is a spanning tree of $\tilde{H'}$, and
hence $H'$ is connected.
Also, by construction, $H'$ is simple, each vertex in \(V_3(H')\) has exactly
one neighbor in \(V_2(H')\), and \(V_2(H')\) is a stable set.
It remains to prove that every cycle in \(H'\) is a separating cycle.
Again, if a cycle \(C' \subset H'\) contains a vertex in \(V_2(H')\) or has a
\(2\)-chord, then \(C'\) is a separating cycle.
Thus, we can assume that \(V(C') \subseteq V_3(H')\) and that \(C'\) has no
\(2\)-chords.
By the construction of $H'$, $C'$ must be a basic cycle of type~1 in $H$, and so
by Claim~\ref{claim:type1cutvertex}, $C'$ is a cut vertex of $\tilde{H}$.
Now we show that $C'$ is a cut vertex in $\tilde{H'}$.
Let $H_1$ be the component in $\tilde{H} \delv C'$ containing the vertex $C$.
Note that, if $C_{i} \in V(\tilde{H} \delv C')$ for $i \in \{1, 2\}$, then $C_{i}$
belongs to $H_{1}$.
Since $\tilde{H'}$ is isomorphic either to $\tilde{H} \delv C$ or to
$\tilde{H} \delv C \adde C_{1}C_{2}$, to show that $C'$ is a cut vertex in
$\tilde{H'}$ it is sufficient to show that $C'$ has neighbor in $V(H_{1} \delv C)$
in the graph $\tilde{H'}$.
If $C'C \notin E(\tilde{H})$, then clearly $C'$ has neighbor in $V(H_{1} \delv C)$.
Thus, we may assume that $C'C \in E(\tilde{H})$.
Now, note that $C' \in \{C_{1}, C_{2}\}$ and $C_{1} \neq C_{2}$, otherwise, by
the construction of $H'$, the cycle $C'$ would be a basic cycle of type~2.
Suppose, without loss of generality, that $C' = C_{1}$.
Therefore, by the construction of $H'$, the edge $C'C_{2} \in E(\tilde{H'})$.
Hence $C'$ has a neighbor in $V(H_{1} \delv C)$ in the graph $\tilde{H'}$, which
implies that $C'$ is a cut vertex in $\tilde{H'}$ and, consequently, that $C'$
is a separating cycle in $H'$.

By induction hypothesis, the graph \(H'\) admits a \(2\)-decomposition
\((F', M')\) satisfying~\eqref{item1} and~\eqref{item2} with respect to the
simple collection of connectors $\mathcal{C}^* \setminus \{y_1\}$.
In what follows, we obtain from \((F',M')\) a \(2\)-decomposition \((F,M)\) of
\(G\) as desired (see Fig.~\ref{fig:decomp2}).

\begin{figure}
\centering
\begin{subfigure}{.5\textwidth}
    \centering
\begin{tikzpicture}[scale = .6]
	
	\node (center1) []	at (0,0) {};
	\node (center2) []	at ($(center1) + (36:4)$) {};
	\node (center3) []	at ($(center1) + (-36:4)$) {};
	\node (center4) []	at ($(center1) + (-72-36:4)$) {};

\foreach \j in {1,...,5}
		{
			\node(1 \j) [black vertex] at ($(center1)+(\j*72+36:1)$) {};
			\node(aux) [black vertex] at ($(center1)+(\j*72+72+36:1)$) {};
			\ifthenelse{\NOT \j=5}{
				\draw[] (1 \j) to [bend right=18] (aux);
			}{
				\draw[color=red] (1 \j) to [bend right=18] (aux);
			}
		}
		
\foreach \j in {1,...,4}
		{
			\node(2 \j) [black vertex] at ($(center2)+(\j*90+36:1)$) {};
			\node(aux) [black vertex] at ($(center2)+(\j*90+90+36:1)$) {};
			\ifthenelse{\NOT \j = 4}{
				\draw[] (2 \j) to [bend right=30] (aux);
			}{
				\draw[color=red] (2 \j) to [bend right=30] (aux);
			}
		}
		
\foreach \j in {1,...,4}
		{
			\node(3 \j) [black vertex] at ($(center3)+(\j*90-36:1)$) {};
			\node(aux) [black vertex] at ($(center3)+(\j*90+90-36:1)$) {};
			\ifthenelse{\NOT \j = 3}{
				\draw[] (3 \j) to [bend right=30] (aux);
			}{
				\draw[color=red] (3 \j) to [bend right=30] (aux);
			}
		}
		
\foreach \j in {1,...,4}
		{
			\node(4 \j) [black vertex] at ($(center4)+(\j*90-18:1)$) {};
			\node(aux) [black vertex] at ($(center4)+(\j*90+90-18:1)$) {};
			\ifthenelse{\NOT \j = 4}{
				\draw[] (4 \j) to [bend right=30] (aux);
			}{
				\draw[color=red] (4 \j) to [bend right=30] (aux);
			}
		}


	\node(in1) [black vertex] at ($.5*(1 1) + .5*(1 2) + (-45:.4)$) {};
	\draw[] 	(1 1) -- (in1);
	\draw[color=red]		(in1) -- (1 2);


	\node(mid12) [black vertex] at ($.5*(1 5)	+.5*(2 2)$) {};
	\draw[] 	(1 5) -- (mid12);
	\draw[color=red] 	(mid12) -- (2 2);

	\node(mid13) [black vertex] at ($.5*(1 4)	+.5*(3 2)$) {};
	\draw[] 	(1 4) -- (mid13);
	\draw[]		(mid13) -- (3 2);
	
	\node(mid14) [black vertex] at ($.5*(1 3)	+.5*(4 1)$) {};
	\draw[color=red] 	(1 3) -- (mid14);
	\draw[]		(mid14) -- (4 1);

	\node () [] at ($(1 3) + (60:.3)$) {$u_k$};
	\node () [] at ($(1 4) + (-90:.4)$) {$x_1$};
	\node () [] at ($(1 5) + (90:.4)$) {$x_2$};
	\node () [] at ($(1 1) + (100:.4)$) {$x_3$};
	\node () [] at ($(1 2) + (180:.4)$) {$z_3$};
	\node () [] at ($(in1) + (-30:.4)$) {$y_3$};
	
	\node () [] at ($(mid12) + (115:.4)$) {$y_2$};
	\node () [] at ($(2 2) + (45:.4)$) {$z_2$};
	
	\node () [] at ($(mid13) + (250:.4)$) {$y_1$};
	\node () [] at ($(3 2) + (-45:.4)$) {$z_1$};

	\node () [] at ($(center1) + (225:1.5)$) {$C$};
	\node () [] at ($(center3) + (-45:1.5)$) {$C_1$};
	\node () [] at ($(center2) + (-45:1.5)$) {$C_2$};
\end{tikzpicture}
    \caption{}
    \label{fig:examplecase2c}
\end{subfigure}%
\begin{subfigure}{.5\textwidth}
    \centering
\begin{tikzpicture}[scale = .6]
	
	\node (center1) []	at (0,0) {};
	\node (center2) []	at ($(center1) + (36:4)$) {};
	\node (center3) []	at ($(center1) + (-36:4)$) {};
	\node (center4) []	at ($(center1) + (-72-36:4)$) {};

\foreach \j in {1,...,5}
		{
			\node(1 \j) [gray vertex] at ($(center1)+(\j*72+36:1)$) {};
			\node(aux) [gray vertex] at ($(center1)+(\j*72+72+36:1)$) {};
			\draw[color=gray] (1 \j) to [bend right=18] (aux);

		}
		
\foreach \j in {1,...,4}
		{
			\node(2 \j) [black vertex] at ($(center2)+(\j*90+36:1)$) {};
			\node(aux) [black vertex] at ($(center2)+(\j*90+90+36:1)$) {};
			\ifthenelse{\NOT \j = 4}{
				\draw[] (2 \j) to [bend right=30] (aux);
			}{
				\draw[color=red] (2 \j) to [bend right=30] (aux);
			}
		}
		
\foreach \j in {1,...,4}
		{
			\node(3 \j) [black vertex] at ($(center3)+(\j*90-36:1)$) {};
			\node(aux) [black vertex] at ($(center3)+(\j*90+90-36:1)$) {};
			\ifthenelse{\NOT \j = 3}{
				\draw[] (3 \j) to [bend right=30] (aux);
			}{
				\draw[color=red] (3 \j) to [bend right=30] (aux);
			}
		}
		
\foreach \j in {1,...,5}
		{
			\node(4 \j) [black vertex] at ($(center4)+(\j*72:1)$) {};
			\node(aux) [black vertex] at ($(center4)+(\j*72+72:1)$) {};
			\ifthenelse{\NOT \j = 4}{
				\draw[] (4 \j) to [bend right=18] (aux);
			}{
				\draw[color=red] (4 \j) to [bend right=18] (aux);
			}
		}


	\node(in1) [gray vertex] at ($.5*(1 1) + .5*(1 2) + (-45:.4)$) {};
	\draw[color=gray] 	(1 1) -- (in1);
	\draw[color=gray]		(in1) -- (1 2);


	\node(mid12) [gray vertex] at ($.5*(1 5)	+.5*(2 2)$) {};
	\draw[color=gray] 	(1 5) -- (mid12);
	\draw[color=gray] 	(mid12) -- (2 2);

	\node(mid13) [gray vertex] at ($.5*(1 4)	+.5*(3 2)$) {};
	\draw[color=gray] 	(1 4) -- (mid13);
	\draw[color=gray]		(mid13) -- (3 2);
	
	\node(mid14) [black vertex] at ($.5*(1 3)	+.5*(4 1)$) {};
	\draw[color=gray] 	(1 3) -- (mid14);
	\draw[color=red]		(mid14) -- (4 1);
	\draw[]		(mid14) -- (4 5);
	
	\node(mid23) [black vertex] at ($.5*(mid12)	+.5*(mid13)$) {};
	\draw[color=red] 	(2 2) -- (mid23);
	\draw[]		(mid23) -- (3 2);
	
	
	\node () [] at ($(2 2) + (45:.4)$) {$z_2$};
	
	\node () [] at ($(3 2) + (-45:.4)$) {$z_1$};

	\node () [] at ($(mid23) + (0:.4)$) {$y$};
	
	\node () [] at ($(center1) + (225:1.5)$) {\color{gray} $C$};
	\node () [] at ($(center3) + (-45:1.5)$) {$C_1$};
	\node () [] at ($(center2) + (-45:1.5)$) {$C_2$};
\end{tikzpicture}
    \caption{}
    \label{fig:examplecase2d}
\end{subfigure}%

\caption{Reduction from a graph $H$ \eqref{fig:examplecase2c} to the graph
    $H'$ \eqref{fig:examplecase2d} in Case 2.
    In both figures, for a \(2\)-decomposition \((F,M)\), the edges in
    \(F\) (resp.\ \(M\)) are colored black (resp.\ red).
    In~\eqref{fig:examplecase2d}, we use gray to indicate the elements from
    $H$ that we removed to create $H'$.}%
\label{fig:decomp2}
\end{figure}
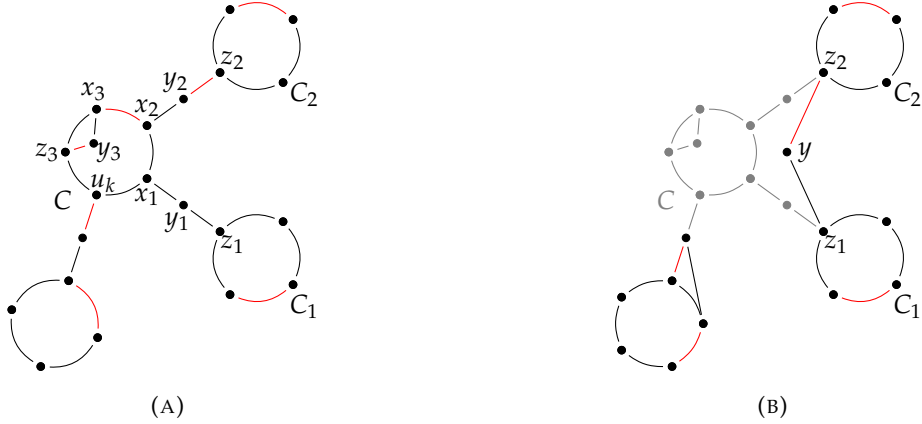

\begin{itemize}
  \item[(a)] We put $x_2x_3$ in $M$ and all the edges of $C \dele x_2x_3$ in $F$.
        We put $y_3x_3$ in $F$, $y_3z_3$ in $M$ and, for each
        $y_i \in J \setminus \{y_3\}$, we put $x_iy_i$ and $y_iz_i$ in distinct
        elements of $(F, M)$.
  \item[(b)] We put $x_1y_1$ and $x_2y_2$ in $F$.
        In addition, for each $y_i \in I \setminus \{y_1,y_2\}$, we put 
        $x_iy_i$ in $M$ and  $y_iz_i$ in $F$.
  \item[(c)] We put each edge
     $$ e \in E(G) \setminus \left( E(C) \cup \big\{z_iz^+_i \colon  y_i \in I \setminus \{y_1,y_2\}\big\} \cup \{y_iz_i, y_ix_i \colon  y_i \in I \cup J\} \right)  $$
        in $F$ if $e \in E(F')$.
        Otherwise, we put $e$ in $M$.
  \item[(d)] For each $y_i \in I \setminus \{y_1,y_2\}$, we put $z_iz^+_i$ in $M$
        if $z'_iz^+_i \in M'$.
        Otherwise, we put $z_iz^+_i$ in $F$.
  \item[(e)] Finally, for each $y_i \in \{y_1,y_2\}$, we put $y_iz_i$ in $F$ if
        $yz_i \in E(F')$.
        Otherwise, we put $y_iz_i$ in $M$.
\end{itemize}

We now show that \((F, M)\) is a \(2\)-decomposition of $H$ as desired.
It is straightforward that each edge of $E(H)$ is either in $F$ or $M$.
Analogously to \textbf{Case 1}, the following claim arises (the same proof
applies).

\begin{claim}\label{claim:notinM2}
  Edge $z_iz'_i \in E(F')$ for all $y_i\in I \setminus \{y_1,y_2\}$.
\end{claim}

From the assignments in steps (a)-(e), since item~(\ref{item1}) holds for
$H'$, and by Claim~\ref{claim:notinM2}, it is clear that $M$ is a matching.
We now check that $F$ is a forest.
Due to the assignments in steps (a)-(d), it is clear that $F$ is a forest in
the graph $G  \deleset  \{y_1z_1, y_2z_2\}$.
The assignments in step~(e) ensure that edges $y_1z_1, y_2z_2$ are assigned
to $F$ or $M$ without creating cycles in $F$: both $y_1z_1, y_2z_2$ are put in
$F$ if $yz_1, yz_2$ are in $F'$, and one of them $y_1z_1$, or $y_2z_2$ is in $M$
whenever its copy in $H'$, namely $yz_1$ or $yz_2$, is in $M'$.
Note that item~(\ref{item1}) holds due to the assignments in steps (a) and
(c), and that item~(\ref{item2}) follows from the assignments in steps
(a)-(c) and Claim~\ref{claim:notinM2}.
This finishes the proof of Case 2, and concludes the proof of
Proposition~\ref{prop:main-stronger}.
\end{proof}

\section{Structural properties of a minimum counterexample}%
\label{sec:counterexample}

In this section, we study the structure of a minimum
counterexample to Conjecture~\ref{con:2dc}.
In fact, in order to avoid technical issues, we work with the following conjecture, which is equivalent to Conjecture~\ref{con:2dc},
as we show in Proposition~\ref{prop:equivalency-2dc-2dc-strong}. 

\begin{conjecture}\label{con:2dc-strong}
  Every graph in $\subcubicIIDCfamily$ can be decomposed into a spanning tree and a matching.
\end{conjecture}

Conjecture~\ref{con:2dc-strong} is more convenient in a proof by induction or minimum counterexample, since it allows vertices with 
degree $1$ without leaving the studied class.
However, the properties of a minimum counterexample
to Conjecture~\ref{con:2dc-strong} can be easily transferred to a minimum counterexample to
Conjecture~\ref{con:2dc} (see \emph{Proof of Theorem~\ref{thm:smallcount2dc}} below).

Let \(\fC\) be the set of all counterexamples for
Conjecture~\ref{con:2dc-strong}, i.e., 
the set of graphs in \(\sepcycle{1}{3}\) that cannot be decomposed into a forest and a matching.
Recall from~\eqref{eq:varphi} that \(\varphi(G) = |V(G)| + |E(G)|\).
Let \(\varphi^* = \min \{\varphi(G) \: G \in \fC\}\), and let
\(\fM = \{ G \in \fC \: \varphi(G) = \varphi^*\}\).
So $\fM$ is the set of all the minimum counterexamples for Conjecture~\ref{con:2dc-strong} according to the function $\varphi$.

First, we show that Conjectures~\ref{con:2dc-strong} and~\ref{con:2dc} are equivalent.
\begin{proposition}\label{prop:equivalency-2dc-2dc-strong}
  Conjecture~\ref{con:2dc-strong} is equivalent to Conjecture~\ref{con:2dc}.
\end{proposition}
\begin{proof}
  Clearly, since \(\sepcycle{2}{3}\subseteq\sepcycle{1}{3}\), Conjecture~\ref{con:2dc} holds if Conjecture~\ref{con:2dc-strong} holds.
  Thus, suppose that Conjecture~\ref{con:2dc} holds and, towards a contradiction, suppose that Conjecture~\ref{con:2dc-strong} does not.
  Let $G \in \fM$ and suppose that \(G\) contains a vertex \(u\) with degree~\(1\).
  Let \(G' = G \delv u\), and note that \(G'\) is a connected subcubic graph in
  \(\subcubicIIDCfamily\) for which \(\varphi(G') < \varphi(G)\).
  Thus, by the minimality of \(G\), the graph \(G'\) contains a 2-decomposition
  \((T', M')\).
  Let \(v\) be the only neighbor of \(u\) in \(G\). 
  Since \(T'\) is
  spanning, \(v \in V(T')\).
  Thus \(T = T'  \adde  uv\) is a spanning tree of \(G\) and \((T, M')\) is a
  2-decomposition of \(G\), a contradiction.
  Thus, we may assume that $V_1(G) = \emptyset$, and so $G$ is a graph in $\iiDCfamily$ that does not admit a 2-decomposition, a contradiction.
\end{proof}

The main theorem of this section is the following.

\begin{theorem}[Structure of a Minimum Counterexample to Conjecture~\ref{con:2dc-strong}]\label{theo:mainsmallestcount}
    Every graph $G \in \fM$ satisfies the following properties:
    \begin{enumerate}
        \item $G$ is simple and 2-edge-connected;
        \item the girth of $G$ is at least 5;
        \item the distance between vertices in \(V_2(G)\) is at least 3;
        \item $G \delvset V_2(G)$ is connected; 
        \item There is a cycle $C \subseteq G$ for which \(V(C) \subseteq V_3(G)\).
    \end{enumerate}
\end{theorem}

Before proving Theorem~\ref{theo:mainsmallestcount}, we show that it implies Theorem~\ref{thm:smallcount2dc}.

\begin{proof}[Proof of Theorem~\ref{thm:smallcount2dc}]
  Suppose $G$ is a minimum counterexample to the 2DC and does not satisfy the
  properties stated in Theorem~\ref{thm:smallcount2dc}. 
  Hence, $G$ is in $\fC$ but is not in $\fM$. 
  Therefore, there is a graph $G' \in \fM$ such that $\varphi(G')<\varphi(G)$. 
  Due to~Claim~\ref{lem:core-no-cuting-edge-vertex} (see below), 
  $G'$ has minimum degree at least 2, and hence $G'$ is a minimum counterexample to the 2DC, a contradiction to the choice of $G$.
\end{proof}

Let $G$ be a graph and let $S \subseteq V(G)$.
The \emph{shrink} of $S$ in the graph $G$, denoted by $G \shrink S$, is the graph obtained from \(G\) 
by identifying the vertices of \(S\) and removing the possible loops.
Given a path \(P\), the \emph{length} of \(P\)
is its number of edges.
A \emph{shortest} path joining vertices \(u\) and \(v\)
is a path joining \(u\) and \(v\) with minimum length
among all such paths;
and the \emph{distance} between $u$ and $v$, denoted by ${\rm dist}_G(u, v)$,
is the length of such a shortest path -- again,  when $G$ is clear from the
context, we may drop the subscript.

\begin{proof}[Proof of Theorem~\ref{theo:mainsmallestcount}]
  \setcounter{claim}{0}
  Let \(G \in \fM\).  First, we show that 
 
  \begin{claim}\label{lem:core-no-cuting-edge-vertex}
    $G$ is 2-edge-connected.
  \end{claim}
  \begin{claimproof}
      Towards a contradiction, suppose that \(G\) contains a cut edge \(e = uv\).
      Let \(G' = G \deleset e\), and let \(G'_u\) and \(G'_v\), respectively, be
      the component of \(G'\) containing the vertex \(u\) and \(v\).
      Clearly, for all \(x \in \{u, v\}\), the graph \(G'_x\) is a connected subcubic
      graph with \(\varphi(G'_x) < \varphi(G)\).
      Moreover, note that \(d_{G'_x}(x) \leq 2\), and hence every cycle containing the
      vertex \(x\) in \(G'_x\) is a separating cycle.
      Thus, it is not hard to check that \(G'_x \in \subcubicIIDCfamily\), and hence, by the
      minimality of \(G\), there are 2-decompositions \((T'_u, M'_u)\) and
      \((T'_v, M'_v)\) of \(G'_u\) and \(G'_v\), respectively.
      Clearly, \((T'_u \cup e \cup T'_v, M'_u \cup M'_v)\) is a 2-decomposition of
      \(G\), a contradiction to the choice of~\(G\).
  \end{claimproof}
  
  Note that since \(G\) has maximum degree at most \(3\), Claim~\ref{lem:core-no-cuting-edge-vertex} implies that \(G\) is \(2\)-connected, and hence has no cut vertex.
  Now, we can prove that $G$ is a simple graph. 
  \begin{claim}\label{lem:core-simple}
    $G$ contains no parallel edges.
  \end{claim}
  \begin{claimproof}
    Towards a contradiction, suppose that \(G\) contains parallel edges.
    Let \(e = uv\) and \(f = uv\) be two parallel edges.
    The graph \(G\) cannot have another parallel edge \(f' = uv\), otherwise the
    \(2\)-cycle \(e \cup f\) would not be a separating cycle.
    If \(d(u) = d(v) = 2\), then \(G\) is a \(2\)-cycle, a contradiction to the
    choice of \(G\).
    Thus, we may assume that \(d(u) = 3\) and that \(u\) has a neighbor \(w\)
    distinct of \(v\).
    Since \(G \in \subcubicIIDCfamily\) and \(e \cup f\) is a cycle, the graph
    \(G \deleset \{e, f\}\) is disconnected, and hence \(uw\) is a cut edge, a
    contradiction to Claim~\ref{lem:core-no-cuting-edge-vertex}.
  \end{claimproof}

  \begin{claim}\label{lem:sc-triangle-free}
      $G$ is triangle-free.
  \end{claim}
  \begin{claimproof}
    Towards a contradiction, suppose that \(G\) contains a triangle \(H = xyzx\).
    By Claim~\ref{lem:core-simple}, \(G\) is simple.
    We may assume that there is a vertex in \(V(H)\), say \(x\), that has
    degree~\(3\) in \(G\), otherwise \(G\) would be a triangle, a contradiction to
    the choice of \(G\).
    If \(d(y) = d(z) = 2\), then \(x\) is a cut vertex, a contradiction to Claim~\ref{lem:core-no-cuting-edge-vertex}.
    Thus, we may assume, without loss of generality, that \(d(y) = 3\).
    The remaining proof is divided in two cases, depending on whether
    \(d(z) = 2 \) or \(d(z) = 3\).

    First, suppose that \(d(z) = 2\).
    Let \(G' = G   \shrink  V(H)\) and name \(u\) the vertex yielded by the shrink of
    \(V(H)\).
    Note that \(G'\) is a connected subcubic graph for which
    \(\varphi(G') < \varphi(G)\).
    Moreover, note that \(d_{G'}(u) = 2\), and hence every cycle in \(G'\)
    containing \(u\) is separating.
    Thus, it is not hard to see that \(G' \in \subcubicIIDCfamily\), and hence, by the minimality
    of \(G\), it follows that there exists a \(2\)-decomposition \((T', M')\) of
    \(G'\).
    For \(v \in \{x, y\}\), let \(e_v\) be the edge in \(G\) incident to \(v\)
    that is not in \(E(H)\), and let \(e'_v\) be the edge in \(G'\) yielded
    from \(e_v\) by the shrinking of \(V(H)\).
    If \(\{e'_x, e'_y\} \subseteq E(T')\), then let
    \(T = T' \deleset \{e'_x, e'_y\} \addeset \{e_x, e_y, xy, yz\}\), and hence
    \((T, M' \cup \{xz\})\) is a \(2\)-decomposition of \(G\), a contradiction to
    the choice of \(G\).
    Thus, we may assume, without loss of generality, that \(e'_x \in M'\), and
    hence \(e'_y \in E(T')\).
    Let \(T = T' \deleset \{e'_y\} \addeset \{e_y,xy,xz\}\) and
    \(M = (M' \setminus \{e'_x\}) \cup \{e_x, yz\}\), and hence \((T, M)\) is a
    2-decomposition of \(G\), a contradiction to the choice of \(G\).

    Now, suppose that \(d(z) = 3\).
    Since \(G \in \subcubicIIDCfamily\), the graph \(G' = G  \deleset  E(H)\) is disconnected.
    It not hard to see that there exists a component in \(G'\) that
    contains only one vertex, say \(x\), in the triangle \(H\).
    Thus, \(x\) is a cut vertex of \(G\), a contradiction to Claim~\ref{lem:core-no-cuting-edge-vertex}.
    \end{claimproof}
  
  \begin{claim}\label{lem:core-no-22-edges}
     If \(u, v \in V_2(G)\), then \({\rm dist}(u, v) \geq 2\).
  \end{claim}
  \begin{claimproof}
    
    Towards a contradiction, suppose that \(e = uv\) is an edge of $G$.
    By Claim~\ref{lem:core-simple}, we have \(|N(u)| \geq 2\). 
    Thus, we may assume that \(N(u) = \{v, x\}\) and \(N(v) = \{u, y\}\) (possibly $x = y$).
    Let \(G' = G \shrink \{u, v\}\), and let \(w\) be the vertex yielded by the contraction of
    the edge \(e\).
    Clearly \(G'\) is a connected subcubic graph for which
    \(\varphi(G') < \varphi(G)\).
    Since \(d_{G'}(w) = 2\), every cycle in \(G'\) containing \(w\) is separating.
    Thus, it is not hard to see that \(G' \in \subcubicIIDCfamily\), and by the minimality of
    \(G\), there exists a 2-decomposition \((T', M')\) of \(G'\).
    Let \(T\) be the spanning tree of \(G\) obtained from \(T'\) by replacing the
    edges incident to \(w\) to its corresponding edges in \(G\) (the edges
    incident to either \(u\) or \(v\)), and by adding the edge \(uv\) to \(T\).
    Therefore \((T, M')\) is a 2-decomposition of \(G\), a contradiction to the
    choice of~\(G\).
  \end{claimproof}

  \begin{claim}\label{lem:sc-deg-2-disj-neigh}
      If \(u \in V_2(G)\) and \(N(u) = \{x,y\}\), then \(N(x) \cap N(y) = \{u\}\).
  \end{claim}
  \begin{claimproof}

    By Claim~\ref{lem:core-simple}, we may assume that \(G\) is simple, thus $x \neq y$, and by Claim~\ref{lem:core-no-22-edges}, we may assume that \(d(x) = d(y) = 3\).
    By Claim~\ref{lem:sc-triangle-free}, \(G\) is triangle-free, and hence the
    edge \(xy \notin E(G)\).

    First, we show that \(|N(x) \cap N(y)| < 3\).
    Towards a contradiction, suppose that \(N(x) = N(y) = \{u, v, w\}\).
    Note that \(vw \notin E(G)\), since \(G\) is triangle-free.
    Let \(S = \{u, x, y, v, w\}\), and note that \(G[S]\) is isomorphic to
    \(K_{2, 3}\).
    We may assume, without loss of generality, that \(d(v) = d(w) = 3\), otherwise
    either \(G\) would be isomorphic to \(K_{2, 3}\), a contradiction to the
    choice of \(G\), or \(G\) would have a cut vertex, a contradiction to Claim~\ref{lem:core-no-cuting-edge-vertex}.
    For \(z \in \{v, w\}\), let \(e_z\) be the edge incident to \(z\) not in
    \(E(G[S])\).
    Let \(G' = G   \shrink  S\), and note that \(G'\) is a connected subcubic graph
    for which \(\varphi(G') < \varphi(G)\).
    Let \(u'\) be the vertex in \(G'\) yielded by the shrink of the set \(S\), and
    note that \(u'\) has degree~\(2\) in \(G'\), and hence, every cycle in \(G'\)
    containing the vertex \(u'\) is separating.
    Thus, it is not hard to see that \(G' \in \subcubicIIDCfamily\), and hence, by the
    minimality of \(G\), there exists a 2-decomposition \((T', M')\) of \(G'\).
    Let \(e'_v\) and \(e'_w\) be the edge of \(G'\) yielded from the edge \(e_v\)
    and \(e_w\), respectively, by the shrink of \(S\).
    If \(\{e'_v, e'_w\} \subseteq E(T')\), then let
    \(T = T' \deleset \{e'_v, e'_w\} \addeset \{e_v, e_w, wy, yv, wx, xu\}\), and hence
    \((T, M' \cup \{uy, xv\})\) is a 2-decomposition of \(G\), a contradiction to
    the choice of \(G\).
    Thus, we may assume, without loss of generality, that \(e'_v \in E(T')\) and
    \(e'_w \in M'\).
    Therefore, let
    \(T = T' \deleset \{e'_v, e'_w\} \addeset \{e_v, e_w, vy, yw, wx, xu\}\), and hence
    \(\big(T, (M' \setminus \{e'_w\}) \cup \{e_w, uy, vx\}\big)\) is a 2-decomposition of
    \(G\), a contradiction to the choice of \(G\).

    Now, we show that \(|N(x) \cap N(y)| < 2\).
    Towards a contradiction, suppose that \(N(x) \cap N(y) = \{u, v\}\).
    Let \(S = \{u, x, y, v\}\), and note that \(G[S]\) is isomorphic to \(C_4\).
    For \(z \in \{x, y, v\}\), let \(e_z\) be the edge incident to \(z\) in \(G\)
    not in \(E(G[S])\) (when \(z = v\), it is possible that the edge \(e_z\) does
    not exist, in this case, let \(e_z\) undefined).
    Let \(G' = G   \shrink  S\), and note that \(G'\) is a connected subcubic graph
    for which \(\varphi(G') < \varphi(G)\).
    Let \(u'\) be the vertex in \(G'\) yielded by the shrink of the set \(S\).
    The remaining proof is divided into two cases depending on whether \(u'\) has
    degree~\(2\) or~\(3\) in \(G'\).

    First, suppose that \(u'\) has degree~\(2\) in \(G'\), and hence, \(v\) has
    degree~\(2\) in \(G\).
    Note that every cycle in \(G'\) containing the vertex \(u'\) is separating.
    Thus, it is not hard to see that \(G' \in \subcubicIIDCfamily\), and hence, by the
    minimality of \(G\), there exists a 2-decomposition \((T', M')\) of \(G'\).
    Let \(e'_x, e'_y,\) be the edges of \(G'\) yielded from the edges \(e_x, e_y\),
    respectively, by the shrink of \(S\).
    Let \(T\) be the subgraph of \(G\) obtained from \(T'\) by replacing each of its
    edges of the type \(e'_z\) for \(e_z\), and by adding the edges \(xv\), \(vy\), and
    \(yu\).
    Let \(M\) be the subgraph of \(G\) obtained from \(M'\) by replacing each of its
    edges of the type \(e'_z\) for \(e_z\), and by adding the edge \(ux\).
    It is not hard to check that \((T, M)\) is a 2-decomposition of \(G\), a
    contradiction to the choice of \(G\).

    Now, suppose that \(u'\) has degree~\(3\) in \(G'\), and hence \(v\) has
    degree~\(3\) in \(G\).
    Let \(a, b, c\) be the neighbors of \(x,y,v\), respectively, in \(G\) which is not in  \(S\).
    Since \(G\) is triangle-free, \(a, b, c\) are three distinct vertices.
    Let \(G' = G \delvset S \adde bc\) and note that \(G'\) is a subcubic
    graph for which \(\varphi(G') < \varphi(G)\).
    Moreover, note that \(G'\) is connected, otherwise, since \(b\) and \(c\) are
    in the same component of \(G'\), \(x\) would be a cut vertex of \(G\), a
    contradiction to~Claim~\ref{lem:core-no-cuting-edge-vertex}.

    Now, we claim that \(G' \in \subcubicIIDCfamily\).
    Towards a contradiction, suppose that \(G'\) contains a non-separating cycle
    \(C'\).
    It is not hard to check that \(bc \in E(C')\).
    Let $P = C' \dele bc$.
    Since $C'$ is a non-separating cycle, there is a path $Q_{z, w}$ in $G'  \deleset  E(C') = G \delvset S  \deleset  E(P)$ for every pair of vertices in $z, w \in V(G) \setminus S$.
    Let $C = P \cup cvyb$ and note that $C$ is a cycle of $G$. 
    Now we show that $C$ is a non-separating cycle in $G$.
    Since there is the path $Q_{z, w}$ in $G  \delvset  S  \deleset  E(P) \subseteq G  \deleset  E(C)$ for every pair of vertices $z, w \in G  \delvset  S$, to prove that $G  \deleset  E(C)$ is connected, it is sufficient to show that there is a path in $G  \deleset  E(C)$ from every vertex in $S$ to $c \in V(G) \setminus S$. 
    Since $G'  \deleset  E(C')$ is connected, there is a path $R$ in $G'  \deleset  E(C') = G  \delvset  S  \deleset  E(P)$ joining 
    $a$ and $c$.
    Moreover, note that $G  \deleset  E(C)$ contains a path from $a$ to every vertex in $S$.
    Therefore,  $G  \deleset  E(C)$ is connected, a contradiction to the choice of $G$.

    Therefore, we may assume that \(G \in \subcubicIIDCfamily\), and hence, by the minimality of
    \(G\), there is a 2-de\-com\-po\-si\-ti\-on \((T', M')\) of \(G'\).
    If \(bc \in E(T')\), then let
    \(T = T' \dele \{bc\} \addeset \{cv, vy, yb, ax, xa\}\), and hence
    \((T, M' \cup \{vx, yu\})\) is a 2-decomposition of \(G\), a contradiction to
    the choice of \(G\).
    Thus, we may assume that \(bc \in M'\), and hence let
    \(T = T' \addeset \{ax, xv, vy, yu\}\), and hence
    \(\big(T, (M' \setminus \{e\}) \cup \{cv, by, xu\}\big)\) is a 2-decomposition of
    \(G\), a contradiction to the choice of \(G\).
  \end{claimproof}

  Finally, we can prove that $G$ has girth at least 5 and that the distance between any two vertices in \(V_2(G)\) is at least 3.

  \begin{claim}\label{lem:sc-girth5}
      \(G\) has girth at least~\(5\).
  \end{claim}
  \begin{claimproof}
    
    By~Claim~\ref{lem:core-simple}, \(G\) is simple, and by~Claim~\ref{lem:sc-triangle-free}, the girth of \(G\) is a least~\(4\).
    Towards a contradiction, suppose that there is a cycle \(C\subseteq G\) with
    length~\(4\).
    Note that \(C\) is an induced cycle.
    By~Claim~\ref{lem:sc-deg-2-disj-neigh}, it
    follows that \(d_G(u) = 3\) for every \(u \in V(C)\).
    
    Let \(G' = G  \deleset  E(C)\).
    Since \(G \in \subcubicIIDCfamily\), the graph \(G'\) is disconnected.
    If there is a component \(H\) of \(G'\) such that
    \(V(H) \cap V(C) = \{u\}\), then $u$ is a cut vertex in $G$, a contradiction to~Claim~\ref{lem:core-no-cuting-edge-vertex}.
    Thus, it follows that $G'$ has precisely two components $H_1$ and $H_2$, 
    and \(|V(H_1) \cap V(C)| = |V(H_2) \cap V(C)| = 2\).
    Let \(V(H_1) \cap V(C) = \{u_1, u_2\}\) and $V(H_2) \cap V(C) = \{u_3, u_4\}$.
    For $i \in \{1, 2, 3, 4\}$, let $u_iv_i \in E(G  \deleset  E(C))$.

    Note that there are two possible cases: (i) \({\rm dist}_C(u_1, u_2) = 1\) or (ii) \({\rm dist}_C(u_1, u_2) = 2\).
    In the first case, suppose, without loss of generality, that $C = u_1u_2u_3u_4u_1$ and, in the second,
    that $C = u_1u_3u_2u_4u_1$ (see Figure~\ref{fig:girth5}).
    Note that in both cases the graph $G$ contains the paths $v_1u_1u_4v_4$ and $v_2u_2u_3v_3$.
    
    \begin{figure}[htp]
      \centering
      \begin{subfigure}{.5\textwidth}
        \centering
  \begin{tikzpicture}[scale=0.6]


    \colorlet{setfilling}{gray!20}
    \colorlet{setborder}{gray!20}

    \path[draw,use Hobby shortcut,closed=true,fill=setfilling, draw=setborder]
    (-3.2,3) .. (-2,2.8) .. (-2,0) .. (-1.8,-3) .. (-3,-3) ;

    \path[draw,use Hobby shortcut,closed=true,fill=setfilling, draw=setborder]
    (2,3) .. (3, 3) .. (3, -3) .. (2,-3) .. (2,0) ;

    \foreach \i in {1,...,4} {
      \node[black vertex] (u_\i) at (\i*90+45:1)  {};
    }

    \node[black vertex] (v_2) at ($(u_2)+(-150:2)$)  {};
    \node[black vertex] (v_3) at ($(u_3)+(-30:2)$)  {};
    \node[black vertex] (v_1) at ($(u_1)+(150:2)$)  {};
    \node[black vertex] (v_4) at ($(u_4)+(30:2)$)  {};

    \foreach[evaluate = \i as \ip using {int(\i + 1)}] \i in {1,...,4} {
      \draw[edge] (u_\i) -- (v_\i);}

    \foreach \i[evaluate = \i as \ip using {int(mod(int(\i + 1), 4))}] in {0,...,3}{
      \pgfmathsetmacro\x{int(\i + 1)}
      \pgfmathsetmacro\y{int(\ip + 1)}
      \draw[edge] (u_\x) -- (u_\y);}

    \foreach \i in {u_1, v_1, u_4, v_4} {
      \node[] at ($(\i)+(90:.3)$) {$\i$};
    }

    \foreach \i in {u_3, v_3, u_2, v_2} {
      \node[] at ($(\i)+(-90:.3)$) {$\i$};
    }
  \end{tikzpicture}%
        \caption{}
        \label{fig:girth5a}
      \end{subfigure}%
      \begin{subfigure}{.5\textwidth}
        \centering
       
  \begin{tikzpicture}[scale=0.6]
    \colorlet{setfilling}{gray!20}
    \colorlet{setborder}{gray!20}
    \path[draw,use Hobby shortcut,closed=true,fill=setfilling, draw=setborder]
    (-3.2,3) .. (-2,2.8) .. (-2,0) .. (-1.8,-3) .. (-3,-3) ;
    \path[draw,use Hobby shortcut,closed=true,fill=setfilling, draw=setborder]
    (2,3) .. (3, 3) .. (3, -3) .. (2,-3) .. (2,0) ;
    \foreach \i in {1,...,4} {
      \node[black vertex] (u_\i) at (\i*90+45:1)  {};
    }
    \node[black vertex] (v_2) at ($(u_2)+(-150:2)$)  {};
    \node[black vertex] (v_3) at ($(u_3)+(-30:2)$)  {};
    \node[black vertex] (v_1) at ($(u_1)+(150:2)$)  {};
    \node[black vertex] (v_4) at ($(u_4)+(30:2)$)  {};
    \foreach[evaluate = \i as \ip using {int(\i + 1)}] \i in {1,...,4} {
      \draw[edge] (u_\i) -- (v_\i);}
    \draw[edge] (u_4) -- (u_1);
    \draw[edge] (u_3) -- (u_2);
    \draw[edge] (u_4) -- (u_2);
    \draw[edge] (u_1) -- (u_3);
    \foreach \i in {u_1, v_1, u_4, v_4} {
      \node[] at ($(\i)+(90:.3)$) {$\i$};
    }
    \foreach \i in {u_3, v_3, u_2, v_2} {
      \node[] at ($(\i)+(-90:.3)$) {$\i$};
    }
  \end{tikzpicture}%
        \caption{}
        \label{fig:girth5b}
      \end{subfigure}%
      \caption{}
      \label{fig:girth5}
    \end{figure}
    
    Let $G'' = G  \delvset  V(C) \addeset \{v_1v_4, v_2v_3\}$ and note that $G''$ is a subcubic connected graph.
    To see that every cycle in $G''$ is separating, suppose the contrary and let $C''$ be a non-separating cycle of \(G''\).
    Since $v_1v_4, v_2v_3$ is an edge cut, it follows that either $C''$ contains both of these edges or none of them.
    If it contains both of them, then clearly $C''$ is a separating cycle.
    If it contains none of them, then $C'' \subseteq G$ and it is not hard to check that $C''$ would be a non-separating cycle in $G$, a contradiction to $G \in \subcubicIIDCfamily$.
    Therefore, $G'' \in \subcubicIIDCfamily$.
    Moreover, note that $\varphi(G'') < \varphi(G)$ and hence, by the minimality of $G$, there is a 2-decomposition $(T'', M'')$ of $G''$.

    First, suppose that \(|E(T'') \cap \{v_1v_4, v_2v_3\}| = 1\).
    Suppose, without loss of generality, that \(v_1v_4 \in E(T'')\) and
    \(v_2v_3 \in M''\).
    If $C = u_1u_2u_3u_4u_1$, then let \(T = T'' \dele v_1v_4 \addeset \{u_1v_1, u_1u_2, u_2u_3, u_3u_4, u_4v_4\}\).
    Otherwise, we have $C = u_1u_3u_2u_4u_1$, and let \(T = T'' \dele v_1v_4 \addeset \{u_1v_1, u_1u_3, u_2u_3, u_2u_4, u_4v_4\}\).
    Let \(M = (M'' \setminus \{v_2v_3\}) \cup \{u_1u_4, u_2v_2, u_3v_3,\}\).
    Therefore, \((T, M)\) is a 2-decomposition of \(G\), a contradiction to the
    choice of \(G\).
    
    Now, suppose that \(|E(T'') \cap \{v_1v_4, v_2v_3\}| = 2\).
    Let
  $$T = T''  \deleset  \{v_1v_4, v_2v_3\} \addeset \{u_1v_1, u_1u_4, u_4v_4, u_2v_2, u_2u_3, u_3v_3\}.$$
    If $C = u_1u_2u_3u_4u_1$, then let \(M = M'' \cup \{u_1u_2, u_3u_4 \}\).
    Otherwise, we have case $C = u_1u_3u_2u_4u_1$, and let $M = M'' \cup \{u_1u_3, u_2u_4\}$.
    Hence, \((T, M)\) is a 2-decomposition of \(G\), a contradiction to the choice
    of \(G\).
  \end{claimproof}

  \begin{claim}\label{lem:sc-dist-deg-2}
     If \(u, v \in V_2(G)\), then \({\rm dist}(u, v) \geq 3\).
  \end{claim}
  \begin{claimproof}
    By~Claim~\ref{lem:core-simple}, we may assume that \(G\) is simple.
    Let \(x\) and \(y\) be the two neighbors of \(u\).
    By~Claims~\ref{lem:core-no-cuting-edge-vertex} and~\ref{lem:core-no-22-edges}, \({\rm dist}(u, v) \geq 2\),
    and hence \(d(x) = d(y) = 3\).
    Let \(N(x) = \{u, x_1, x_2\}\) and \(N(y) = \{u, y_1, y_2\}\).
    By~Claim~\ref{lem:sc-deg-2-disj-neigh}, the
    vertices \(x_1, x_2, y_1, y_2\) are all distinct.
    Towards a contradiction, suppose that \(d(u, v) = 2\).
    We may assume, without loss of generality, that \(v = x_1\).
    Let \(G' = (G \delv x)   \shrink  \{x_1, u\}\), and let \(u'\) be the vertex in \(G'\)
    yielded by the shrink of the set $\{x_1, u\}$.
    Note that \(G'\) is a subcubic graph for which \(\varphi(G') < \varphi(G)\).
    Moreover, note that \(G'\) is connected, otherwise \(x\) is a cut vertex of
    \(G\), a contradiction to Claim~\ref{lem:core-no-cuting-edge-vertex}.
    Since \(u'\) has degree~\(2\) in \(G'\), every cycle in \(G'\) containing
    \(u'\) is a separating cycle.
    It is not hard to check that in fact \(G' \in \subcubicIIDCfamily\).
    Therefore, by the minimality of \(G\), there exists a 2-decomposition
    \((T', M')\) of \(G'\).
    Let \(e_{x_1}\) and \(e_u\) be the edges in \(G\) incident to \(x_1\) and
    \(u\), respectively, which is not incident to \(x\).
    Let \(e'_{x_1}\) and \(e'_{u}\) be the edge yielded from \(e_{x_1}\) and
    \(e_u\), respectively, by the shrink of \(\{u, x_1\}\) in \(G\).
    The remaining proof is divided into two cases depending on whether
    \(|E(T') \cap \{e'_{x_1}, e'_u\}|\) is~\(1\) or~\(2\).

    First, suppose that \(|E(T') \cap \{e'_{x_1}, e'_u\}| = 1\).
    Suppose, without loss of generality, that \(e'_{x_1} \in E(T')\) and
    \(e'_u \in M'\).
    Thus, let \(T = T' \deleset \{e'_{x_1}\} \addeset \{e_{x_1}, x_2x, xu\}\), and
    hence \(\big(T, (M' \setminus \{e'_u\}) \cup \{e_u, xx_1\}\big)\) is a
    2-decomposition of \(G\), a contradiction to the choice of \(G\).
    Now, suppose that \(|E(T') \cap \{e'_{x_1}, e'_u\}| = 2\).
    Thus \(d_{T'}(u') = 2\), the graph \(T' \delv u'\) has two components.
    Let \(T' = T'_{x_1} \cup T'_u\) where \(V(T'_{x_1}) \cap V(T'_u) = \{u'\}\),
    \(e'_{x_1} \in E(T'_{x_1})\), and \(e'_u \in E(T'_u)\).
    Suppose, without loss of generality, that \(x_2 \in V(T'_{x_1})\), and let
    \(T = T' \deleset \{e'_{x_1}, e'_u\} \addeset \{e_{x_1}, e_u, x_2x, xu\}\). 
    Hence, \((T, M' \cup \{x_1x\})\) is a 2-decomposition of \(G\), a
    contradiction.
  \end{claimproof}

  In what follows we prove the last two properties, namely that $G - V_2(G)$ is connected and that $G$ has a cycle containin only  vertices with degree~$3$.
  For that, let $V_2(G) = \{v_1, v_2, \ldots, v_\ell\}$ and for each $v_i \in V_2(G)$,
  let $N(v_i) = \{u_i , u'_i\}$. 

  \begin{claim}\label{cla:g-v2-connected}
    $G  \delvset  V_2(G)$ is connected.
  \end{claim}
  \begin{claimproof}
      Let $v_i \in V_2(G)$.
      Since $G$ is triangle-free, $u_iu'_i \notin E(G)$.
      Let $G' = G \delv v_i \adde u_iu'_i$ and note that $G'$ is a subcubic connected graph.
      If every cycle in $G'$ is separating, then, by the minimality of $G$, there is a 2-decomposition $(T', M')$ of $G$.
      If $u_iu'_i \in M'$, then  $(T' \adde u_iv_i, M' \cup \{u_i'v_i\})$ is a 2-decomposition of $G$ and, if $u_iu'_i \in E(T')$, then $(T' \addeset \{u_iv_i, u_i'v_i\}, M')$ is a 2-decomposition of $G$.
      In both cases we have a contradiction to the choice of $G$.
      
      Thus, we may assume that $G'$ contains a non-separating cycle $C$.
      Note that $u_iu_i' \in E(C)$.
      Moreover, all vertices in $C$ have degree $3$ in $G'$, and so in $G$, since $d_{G'}(v) = d_G(v)$ for all $v \in V(C)$. 
      Let $P_i = C \dele u_iu_i'$ and note that $P_i$ is a path in $G \delv v_i$ joining the vertices $u_i$ and $u_i'$ containing only vertices in \(V_3(G)\). 
      Since $v_i$ was arbitrarily chosen, theses properties hold for every $v_i \in V_2(G)$.
    If \(G' = G-V_2(G)\) is disconnected,
    then there is \(v_i\in V_2(G)\) for which  \(u_i\) and \(u'_i\),
    are in different components of \(G'\).
    But since \(V(P_i)\subseteq V_3(G)\),
    the path \(P_i\) is in \(G'\) and joins \(u_i\) and \(u'_i\),
    a contradiction.
  \end{claimproof}

  \begin{claim}
      There is a cycle $C \subseteq G$ containing only degree~$3$ vertices.
  \end{claim}
  \begin{claimproof}
      By Claim~\ref{cla:g-v2-connected}, $G'= G  \delvset  V_2(G)$ is connected.
      Towards a contradiction, suppose that $G'$ is a tree.
      Let $T = G' \addeset \{v_iu_i \: v_i \in V_2(G)\}$ and note that $T$ is a spanning tree of $G$.
      Let $M = \{v_iu'_i \: v_i \in V_2(G)\}$ and let $v_ju'_j \in M$.
      Since $d_G(v_j) = 2$ and $v_ju_j \in E(T)$, there is only one edge in $M$
      incident to $v_j$, namely $v_ju'_j$.
      If there is another edge in $M$, say $v_ku'_k$, incident to $u'_j$, then $u'_j = u'_k$ and ${\rm dist}_G(v_j, v_k) \leq 2$, a contradiction to Claim~\ref{lem:sc-dist-deg-2}.
      Thus, $M$ is a matching in $G$ and $(T, M)$ is a 2-decomposition of $G$, 
      a contradiction.
      Therefore, $G  \delvset  V_2(G)$ is not a tree, and hence contains a cycle $C \subseteq G - V_2(G)$.
      By Claim~\ref{lem:core-no-cuting-edge-vertex}, $\delta(G) \geq 2$, and hence $G - V_2(G) = G[V_3]$.
      Therefore, $d_G(u) = 3$ for all $u \in V(C)$.
  \end{claimproof}
\end{proof}

\section{Concluding Remarks}

We recall that a graph is \emph{claw-free} if it contains no induced copy of
\(K_{1, 3}\).
As mentioned before, the 3DC has been verified for claw-free cubic graphs~\cite{HoLiYu20, AbAhAk20}.
Using an idea similar to the one used by Aboomahigir, Ahanjideh and
Akbari~\cite{AbAhAk20}, one can show  that the 2DC also holds for claw-free graphs.
For completeness, we provide its proof below.

\begin{theorem}\label{th:claw-free-2d}
  If \(G \in \subcubicIIDCfamily\) is a claw-free graph, then
  \(G\) can be decomposed into a spanning tree and a matching.
\end{theorem}
\begin{proof}
  Towards a contradiction, suppose the opposite, and let \(G\) be a
  counterexample with a minimum number of vertices.
  It is easy to check that the result holds for $|V(G)| \leq 2$, 
  thus we may assume that $|V(G)| > 2$.
  
  First, suppose that $G$ contains a cut edge $uv$ and let $G' = G \dele uv$.
  Let $G_u$ and $G_v$ be the components of $G'$ containing the vertices $u$ and $v$, respectively.
  By the minimality of $G$, there are 2-decompositions $(T_u, M_u)$ and $(T_v, M_v)$ of $G_u$ and $G_v$, respectively. 
  Hence, $(T_u \cup T_v \cup \{uv\}, M_u \cup M_v)$ is a 2-decomposition of $G$, a contradiction.
  Thus, we may assume that $G$ is $2$-edge connected, and hence $\delta(G) \geq 2$.
  
  Moreover, $G$ must contain a vertex $v$ with degree $3$, otherwise $G$ would be a cycle and the result follows.
  Let \(N(v) = \{x, y, z\}\).
  Since \(G\) is claw-free, we may assume without loss of generality that
  \(xy \in E(G)\).
  Let \(C\) be the cycle \(xvyx\).
  Since \(G \in \subcubicIIDCfamily\), the graph \(G  \deleset  E(C)\) is disconnected.
  If \(d(x) = d(y) = 3\), then \(G\) has a cut edge, a contradiction.
  Hence we may assume without loss of generality that \(d(x) = 2\).
  If \(yz \in E(G)\), then either \(G\) is the diamond ($K_4$ minus an edge) or \(G\) has a cut edge incident to $z$;
  in both cases, we reach a contradiction.
  
  Thus, we may assume that \(yz \notin E(G)\).
  Note that \(d(y) = 3\) and \(d(z) \geq 2\), otherwise \(vz\) would be a cut edge.
  Let \(G' = G \delvset \{x, v\} \addeset yz\), and note that \(G'\) is a connected subcubic graph.
  
  We claim that $G'$ is claw-free.
  Indeed, if $d_G(z) = 2$, then it is easy to see that $G'$ is claw-free.
  Thus suppose that $d_G(z) = 3$ and let $N_G(z) = \{v, a, b\}$. 
  Note that $a, b \notin \{x, y\}$, because $yz \notin E(G)$ and $d_G(x) = 2$.
  Since $\{z, v, a, b\}$ does not induces a claw and $va, vb \notin E(G)$, it follows that $ab \in E(G)$.
  Thus, \(\{z, y, a, b\}\) does not induces a claw in $G'$, and hence $G'$ is claw-free.

  Now we show that $G' \in \subcubicIIDCfamily$.
  Let \(C'\) be a cycle in \(G'\).
  If \(yz \in C'\), then \(C'\) is a separating cycle, since \(d_{G'} (y) = 2\).
  Otherwise, \(yz \notin E(C')\), and hence, \(C'\) is a cycle in
  \(G \in \subcubicIIDCfamily\).
  In this case, it is easy to check that \(C'\) is separating.
  
  Therefore, \(G'\) is a claw-free, connected, subcubic graph in
  \(\subcubicIIDCfamily\) with fewer vertices than \(G\), and by the
  minimality of \(G\), the graph \(G'\) admits a decomposition into a spanning
  tree \(T'\) and a matching \(M'\).

  Now, we show how to find a decomposition of \(G\) into a spanning tree \(T\)
  and a matching \(M\) from \(T'\) and \(M'\).
  If \(yz \in M'\), then let \(T = T \addeset \{yv, vx\}\) and
  \(M = (M' \setminus \{yz\}) \cup \{vz, xv\}\) (recall that
  \(yz \notin E(G)\)).
  Otherwise, \(yz \in T'\), and hence let \(T = T' \dele yz \addeset \{yv, yx, vz\}\) and
  \(M = M' \cup \{vx\}\).
  Note that in both cases \((T, M)\) is a 2-decomposition of $G$, a contradiction.
\end{proof}

\section{Competing interest and Data Availability statement}

The authors declare that they have no conflict of interest.
Data sharing is not applicable to this article as no new data were created or analyzed in this study.

\bibliographystyle{amsplain}

\end{document}